\documentclass[12pt]{amsart}
\usepackage{amssymb}  
\usepackage{latexsym} 
\usepackage[all]{xy}

\newcommand{\Sha}{\mbox{\wncyr Sh}}
\newcommand{\Z}{{\mathbb Z}}
\newcommand{\Q}{{\mathbb Q}}
\newcommand{\R}{{\mathbb R}}
\newcommand{\PP}{{\mathbb P}}
\newcommand{\CL}{{\mathcal L}}
\newcommand{\tensor}{\otimes}
\newcommand{\charic}{\operatorname{char}}
\newcommand{\Kbar}{\overline{K}}
\newcommand{\Fbar}{\overline{F}}
\newcommand{\Magma}{{\sf MAGMA }}
\newcommand{\AGL}{\operatorname{AGL}}
\newcommand{\GL}{\operatorname{GL}}
\newcommand{\WC}{\operatorname{WC}}
\newcommand{\SL}{\operatorname{SL}}
\newcommand{\PGL}{\operatorname{PGL}}
\newcommand{\PSL}{\operatorname{PSL}}
\newcommand{\Gal}{\operatorname{Gal}}
\newcommand{\Aut}{\operatorname{Aut}}
\newcommand{\End}{\operatorname{End}}
\newcommand{\Hom}{\operatorname{Hom}}
\newcommand{\inv}{\operatorname{inv}}
\newcommand{\Ext}{\operatorname{Ext}}
\newcommand{\Jac}{\operatorname{Jac}}
\newcommand{\Map}{\operatorname{Map}}
\newcommand{\Sym}{\operatorname{Sym}}
\newcommand{\Sel}{\operatorname{Sel}}
\newcommand{\Pic}{\operatorname{Pic}}
\newcommand{\Mat}{\operatorname{Mat}}
\newcommand{\eps}{\varepsilon}
\newcommand{\Gm}{\mathbb{G}_{\text{\rm m}}}
\newcommand{\isom}{\cong}
\newcommand{\inj}{\hookrightarrow}
\newcommand{\ra}{\longrightarrow}

\newcommand{\Div}{\operatorname{Div}}
\newcommand{\op}{{\operatorname{op}}}
\newcommand{\Tr}{\operatorname{Tr}}
\newcommand{\tr}{\operatorname{tr}}
\newcommand{\divv}{\operatorname{div}}
\newcommand{\Abar}{{\overline{A}}}
\newcommand{\Rbar}{{\overline{R}}}
\newcommand{\Ob}{\operatorname{Ob}}
\newcommand{\Br}{\operatorname{Br}}
\newcommand{\alphahat}{\widehat{\alpha}}
\newcommand{\ndv}{\!\nmid\!}

\newfont{\wncyr}{wncyr10 at 12pt}
\newfont{\wncyrten}{wncyr10 at 10pt}

\newenvironment{Proof}{\par\noindent{\sc Proof:}}%
                      {\hspace*{\fill}\nobreak$\Box$\par\medskip}

\newtheorem{Proposition}{Proposition}[section]
\newtheorem{Lemma}[Proposition]{Lemma}
\newtheorem{Corollary}[Proposition]{Corollary}

\theoremstyle{definition}
\newtheorem{Definition}[Proposition]{Definition}
\newtheorem{Remark}[Proposition]{Remark}
\newtheorem{Remarks}[Proposition]{Remarks}

\begin{document}

\date{1st March 2006}
\title[Explicit $n$-descent on elliptic curves]%
{Explicit $n$-descent on elliptic curves\\I. Algebra}

\author{J.E.~Cremona}
\address{School of Mathematical Sciences,
         University of Nottingham, 
         University Park, Nottingham NG7 2RD, UK}
\email{John.Cremona@nottingham.ac.uk}

\author{T.A.~Fisher}
\address{University of Cambridge,
         DPMMS, Centre for Mathematical Sciences,
         Wilberforce Road, Cambridge CB3 0WB, UK}
\email{T.A.Fisher@dpmms.cam.ac.uk}

\author{C. O'Neil}
\address{Barnard College, Columbia University,
Department of Mathematics,
2990 Broadway MC 4418,
New York, NY 10027-6902,USA}
\email{oneil@math.columbia.edu}

\author{D. Simon}
\address{Universit\'e de Caen, Campus II - Boulevard Mar\'echal Juin,
BP 5186--14032, Caen, France}
\email{Denis.Simon@math.unicaen.fr}

\author{M. Stoll}
\address{School of Engineering and Science, 
International University Bremen,
P.O. Box 750561, 28725 Bremen, Germany}
\email{M.Stoll@iu-bremen.de}

\begin{abstract} 
This is the first in a series of papers in which 
we study the $n$-Selmer group of an elliptic curve,
with the aim of representing its elements as 
genus one normal curves of degree $n$.
The methods we describe are practical in the case $n=3$ for elliptic
curves over the rationals, and have been implemented in \Magma.
\end{abstract}

\maketitle

\section*{Introduction}

Descent on an elliptic curve $E$, defined over a number field $K$,
is a method for obtaining information about
both the Mordell-Weil group $E(K)$ 
and the Tate-Shafarevich group $\Sha(E/K)$. 
Indeed for each integer $n\ge2$ there is an exact sequence
$$ 0 \to E(K)/nE(K) \to \Sel^{(n)}(E/K) \to \Sha(E/K)[n] \to 0 $$
where $\Sel^{(n)}(E/K)$ is the $n$-Selmer group. 

This is the first in a series of papers in which 
we study the $n$-Selmer group with the aim of representing its 
elements as genus one normal curves $C \subset \PP^{n-1}$ (when $n \ge 3$).
Having this representation allows searching for rational
points on $C$ (which in turn gives points in $E(K)$, since
$C$ may be seen as an $n$-covering of $E$) and is a first step towards
doing higher descents. A further application is to the study of 
explicit counter-examples to the Hasse Principle.

In this introduction we discuss our approach to the problem
and set out the goals of our work. Following some historical remarks, 
we will outline the contents of this first paper, and then briefly that 
of the remaining papers in the series.

The method of descent, for explicitly determining the solutions of
Diophantine Equations, has a long and distinguished history going back
(at least) to Fermat.  As a tool for the determination of the
Mordell-Weil groups of elliptic curves over number fields, descent has
been used since the very first applications of computing to number
theory.  Until the 1990s, the only methods which had been implemented
for general elliptic curves were based on $2$-descent and could be
applied only to elliptic curves defined over $\Q$, though of course
individual examples had been worked out over other fields.  The advent
of higher-level computer algebra software and the development of
efficient algorithms for handling the arithmetic of more general
number fields has meant that $2$-descent can now be carried out over
general number fields (of moderate degree and discriminant, for
practical reasons).  We may cite both Simon's {\tt gp} program {\tt
ell.gp} (see \cite{Simon_ell_paper}), and the {\Magma} package written
primarily by Bruin, as examples of this which are widely used.

The situation regarding so-called higher descents, meaning $n$-descent
for $n>2$, has until now been much more fragmentary and less
satisfactory.  Some $3$-descents (for twists of Fermat curves) and
certain other descents via isogeny have been studied systematically,
but these apply only to special families and not to general elliptic
curves (at least, not without an extension of the ground field, which
introduces further complexities and complicates implementation
significantly).  Higher $2$-power descents (also known as second and
third $2$-descents) have been studied and implemented by Siksek
\cite{MSS} and Womack \cite{WomackThesis} for $4$-descent, and by
Stamminger \cite{StammingerThesis} for $8$-descent.

In the case of $2$-descent, the map $C\to\PP^1$ is a double cover
rather than an embedding, and the elements of $\Sel^{(2)}(E/K)$
are represented as curves of the form $Y^2=g(X)$ where $g$ is a
quartic.  Our goal is to be equally explicit for all $n>2$. 
The methods we present are fully worked out for all odd prime~$n$, 
and have been implemented in \Magma \cite{Magma} in the case $n=3$ 
for elliptic curves over~$\Q$. This implementation will be included in
\Magma version 2.13, to be released later this year.

To avoid making assumptions about the Galois module structure of
$E[n]$, we work with the \'etale algebra $R$ of $E[n]$.
This is a $K$-algebra of dimension $n^2$, 
explicitly realised as a product of number fields. 
The starting point for our work is the paper of Schaefer and
Stoll \cite{SchaeferStoll}, which improves on earlier methods in \cite{DSS}.
They show that if $n$ is prime then a certain group homomorphism 
$$ w_1 : H^1(K,E[n]) \to R^\times / (R^\times)^n $$
is injective, and determine its image. This is the basis of
an algorithm for computing $\Sel^{(n)}(E/K)$ as a subgroup
of $R^\times / (R^\times)^n$. (In fact they assume that $n$ is 
odd, since the case $n=2$ was already well known.)
The algorithm requires knowledge of the class group and unit 
group of each constituent field of $R$.

In \S\ref{etal} we replace $w_1$ by a group
homomorphism 
$$ w_2 : H^1(K,E[n]) \to (R \otimes R)^\times / \partial R^\times $$
where $\partial : R^\times \to (R \otimes R)^\times$ is a certain
map. We show that $w_2$ is injective for all $n \ge 2$ and 
determine its image. If $n$ is prime then it is possible to convert
the subgroup of $R^\times / (R^\times)^n$  computed by 
Schaefer and Stoll to a subgroup of 
$(R \otimes R)^\times / \partial R^\times$.
Since $R \otimes R$ is the \'etale algebra of $E[n] \times E[n]$,
it is again a product of number fields.
In principle one could compute the $n$-Selmer group directly
as a subgroup of $(R \otimes R)^\times / \partial R^\times$,
but this would require knowledge of the class group and unit
group of each constituent field of $R \otimes R$,
and in general these fields have larger degree than those appearing
in $R$. 

Our goal is now the following. We must convert elements of 
$\Sel^{(n)}(E/K)$ represented algebraically by 
$\rho \in (R \otimes R)^\times$,
to elements of $\Sel^{(n)}(E/K)$ represented geometrically by 
(equations for) genus one normal curves $C \subset \PP^{n-1}$. 
We present three algorithms for performing this conversion, the
second and third of which apply for arbitrary $n \ge 2$. 
In particular we have no need to assume that $n$ is prime.
In the case $n=2$ our third algorithm reduces to the classical
number field method for 2-descent. Nevertheless we 
assume for ease of exposition that $n \ge 3$. 

We give a brief description of each algorithm.

\medskip

\paragraph{\bf The Hesse pencil method} We determine $n \times n$ matrices
(with entries in an extension of $K$) that represent the action
of $E[n]$ on $C \subset \PP^{n-1}$. At least in the case 
$n=3$ it is then practical to recover an equation for $C$.

\medskip

\paragraph{\bf The flex algebra method} We embed $E \subset \PP^{n-1}$ via
the complete linear system $|n.0|$ where $0$ is the identity on
$E$. We then determine a change of co-ordinates (defined over 
an extension of $K$) that takes $E \subset \PP^{n-1}$ to
$C \subset \PP^{n-1}$. We use this change of co-ordinates to 
compute equations for $C$.

\medskip

\paragraph{\bf The Segre embedding method} We determine equations for $C$
as a curve of degree $n^2$ in the rank 1 locus of $\PP(\Mat_n)$. 
Thus $C$ lies in the image of the Segre embedding
$$ \PP^{n-1} \times (\PP^{n-1})^\vee \to \PP(\Mat_n). $$
We pull back to $\PP^{n-1} \times (\PP^{n-1})^\vee$ and then project
onto the first factor to get $C \subset \PP^{n-1}$.

\medskip

It is important to realise that not every element of $H^1(K,E[n])$ 
can be represented by a genus one normal curve $C \subset \PP^{n-1}$.
Those elements of $H^1(K,E[n])$ that can be represented in this
way form the ``kernel'' of the obstruction map
$$ \Ob : H^1(K,E[n]) \to \Br(K) $$
where $\Br(K)$ is the Brauer group of $K$. The reader is warned that
in general the obstruction map is not a group homomorphism, and its
kernel is not a group.

In fact each of our algorithms works over an arbitrary field $K$
(assumed perfect and with $\charic(K) \ndv n$) provided we make
the following hypotheses:

\begin{itemize}

\item We start with an element $\rho \in (R \otimes R)^\times$ 
that represents an element of $H^1(K,E[n])$ with trivial obstruction. 

\item We have access to a ``Black Box'' that, given structure
constants for a $K$-algebra known to be isomorphic to $\Mat_n(K)$,
finds such an isomorphism explicitly.

\end{itemize}

Returning to the case $K$ a number field, it follows by the
commutativity of the diagram
$$ \xymatrix@C=3em{ H^1(K,E[n]) \ar[r]^-{\Ob} \ar[d] & \Br(K) \ar[d] \\
\prod_v H^1(K_v,E[n]) \ar[r]^-{\prod_v \Ob_v} & \prod_v \Br(K_v) } $$
and the injectivity of the right hand map, that an element
of $H^1(K,E[n])$ has trivial obstruction if and only if it has
trivial obstruction everywhere locally. We deduce
that $\Sel^{(n)}(E/K)$ is contained in the kernel of the
obstruction map. In our algorithms it then becomes necessary to
use an explicit version of the local-to-global principle for the
Brauer group. This role is played by the Black Box. 
An essentially equivalent problem is that of finding a 
rational point on a Brauer-Severi variety of dimension $n-1$. 
In the case $n=2$ this means finding a rational point on a conic, 
a task which we recognise as one of the steps in the classical 
number field method for $2$-descent, cf. \cite{CaL}, \S15.  

\bigskip

We present our work in a series of papers, of which this is 
the first. We briefly summarise the contents of each.

\medskip

\paragraph{\bf Paper I: Algebra} We work over a perfect field $K$
with $\charic(K) \ndv n$. We give a list of interpretations
of the Galois cohomology group $H^1(K,E[n])$, and explore the
relationships between them.
Then we go through several different descriptions of the obstruction map
and check that they are all equivalent. In \S\ref{etal} we 
introduce the \'etale algebra $R$ of $E[n]$ and define  
the maps $w_1$ and $w_2$. In 
\S\ref{algebras} we explain how the element $\rho \in (R \otimes R)^\times$
may be used to construct certain $K$-algebras. We end by outlining
each of our three algorithms, assuming in each case the existence
of a suitable Black Box. 

\medskip

\paragraph{\bf Paper II: Geometry}
We give further details of the Segre embedding method.
We represent elements of $H^1(K,E[n])$ by 
Brauer-Severi diagrams $[C \to S]$. The obstruction in $\Br(K)$ is 
represented both by the Brauer-Severi variety $S$ (of dimension $n-1$)
and by the obstruction algebra $A$ (a central simple algebra over $K$ of 
dimension $n^2$). We specify certain embeddings 
   $$ C \to S \times S^\vee \to \PP(A) $$
where $S^\vee$ is the dual of $S$. Then starting from
$\rho \in (R \otimes R)^\times$ we explain how to 
write down both structure constants for $A$, and 
equations for $C$ as a curve in $\PP(A)$. If the obstruction is 
trivial then there are isomorphisms $S \isom \PP^{n-1}$ and 
$A \isom \Mat_n(K)$. We recover $C \subset \PP^{n-1}$ by pulling 
back the image of $C$ in $\PP(\Mat_n)$ to $\PP^{n-1} \times (\PP^{n-1})^\vee$ 
and then projecting onto the first factor.

\medskip

\paragraph{\bf Paper III: Algorithms}
We work over a number field $K$ and assume that $n$ is prime.
We briefly recall how to compute 
$\Sel^{(n)}(E/K)$ first as a subgroup of $R^ \times/ (R^\times)^n$ 
and then as a subgroup of $(R  \otimes R)^\times / \partial R^\times$.
Then we give further practical details of our algorithms
for converting $\rho \in (R \otimes R)^\times$ to $C \subset \PP^{n-1}$,
concentrating on the Segre embedding method in the case $n=3$.
We also describe the methods we use for the Black Box,
including one that is practical in the case $K=\Q$ and $n=3$.
We illustrate our work with numerical examples.

\section{Interpretations of $H^1(K,E[n])$}
\label{sect1} 

Let $E$ be an elliptic curve defined over a perfect field $K$. We write 
$G_K=\Gal(\Kbar/K)$, where $\Kbar$ is the algebraic closure of $K$, and
$H^i(K,-)$ for the Galois cohomology group (or set) 
$H^i(G_K,-)$.
Let $n \ge 2$ be an integer with $\charic(K) \ndv n$. 
Taking Galois cohomology of the exact sequence 
$$ 0 \ra E[n] \ra E \stackrel{n}{\ra} E \ra 0 $$
we obtain the Kummer exact sequence
$$ \ldots \to E(K) \to H^1(K,E[n]) \to H^1(K,E) \to \ldots. $$
We discuss how to represent elements of these groups.
For the group on the left this is straight forward:
we fix a Weierstrass equation for $E$ and specify 
points in $E(K)$ by writing down their co-ordinates.

To interpret the other two groups we call
upon the general principle that the twists of 
an object $X$ (defined over $K$) are parametrised by $H^1(K,\Aut(X))$
where $\Aut(X)$ is the automorphism group of $X$.
More precisely if  
$Y$ is another object defined over $K$, and $\phi : Y \to X$ 
is an isomorphism defined over $\Kbar$, then associating to $Y$
the cocycle $\xi_\sigma = \sigma(\phi) \phi^{-1}$, determines
an injective map from the $K$-isomorphism classes of twists 
to $H^1(K,\Aut(X))$. We claim that in each of our applications this map is 
also surjective. Indeed if $X$ is a quasi-projective $K$-variety,
then the surjectivity follows by Galois descent: see 
\cite{SerreAGCF}, Chap.~V, Cor.~2 to Prop.~12.
In general $X$ will be a quasi-projective
$K$-variety $X_0$ equipped with certain ``additional structure''.
Thus $\Aut(X) \subset \Aut(X_0)$. To construct the twist of 
$X$ by $\xi \in H^1(K,\Aut(X))$ we first take the twist $Y_0$ of $X_0$
by $\xi$ and then use the isomorphism $\phi : Y_0 \to X_0$
with $\sigma(\phi) \phi^{-1}= \xi_\sigma$ to transfer the additional
structure on $X_0$ to $Y_0$. A routine calculation shows that
the additional structure on $Y_0$ is Galois invariant, and so defined
over $K$. We recall that if $\Aut(X)$ is abelian then $H^1(K,\Aut(X))$ 
is an abelian group;  otherwise, it is a pointed set with 
the class of~$X$ as its distinguished element.

The prototype example is that of a torsor or 
principal homogeneous space under $E$.
In \S\ref{Entorsors} we will also consider torsors under $E[n]$.

\begin{Definition}
(i) A torsor 
under $E$ is a pair $(C,\mu)$
where $C$ is a smooth projective curve of genus one 
(defined over $K$) and $\mu:E \times C \to C$ 
is a morphism (defined over $K$) 
that induces a simple transitive action on $\Kbar$-points. \\
(ii) An isomorphism of torsors $(C_1,\mu_1) \isom (C_2,\mu_2)$ 
is an isomorphism of curves $C_1 \isom C_2$ that respects the action of $E$.
\end{Definition}

The trivial torsor is $(E,+)$ where $+:E \times E \to E$ is the
group law. 

\begin{Lemma}
\label{torsorsaretwists}
Every torsor under $E$ is a twist of $(E,+)$.
\end{Lemma}
\begin{Proof} More generally we note that if $(C,\mu)$ is a 
torsor under $E$ and $P_0 \in C(\Kbar)$ then $(E,+) \isom (C,\mu) \, ; \,
P \mapsto \mu(P,P_0)$ is an isomorphism of torsors defined over $K(P_0)$.
\end{Proof}

\begin{Lemma} 
\label{auttorsor}
$\Aut(E,+) \isom E$.
\end{Lemma}
\begin{Proof} The automorphisms of $E$, as a torsor under itself, are
the automorphisms of $E$, as a curve, that commute with all translation maps. 
The only such automorphisms are the translation maps themselves.
\end{Proof}

By the twisting principle we obtain

\begin{Proposition} 
\label{prototype}
The torsors under $E$, viewed as twists of $(E,+)$,
are parametrised up to isomorphism by $H^1(K,E)$. 
\end{Proposition}

With this geometric interpretation, the group $H^1(K,E)$ is called
the Weil-Ch\^atelet group $\WC(E/K)$.

\begin{Remark} By a standard abuse of notation we refer to a torsor
as $C$ rather than $(C,\mu)$. The only ambiguity in the choice of 
$\mu$ comes from the automorphisms of~$E$ as an elliptic curve. 
If $j(E) \not= 0,1728$ these are just $\{\pm 1\}$, so $C$ will have at most two
possible structures of torsor under~$E$. The two structures coincide when 
either one has order dividing~$2$ in~$H^1(K, E)$.
\end{Remark}

We have prepared the following list of interpretations of the
group $H^1(K,E[n])$. Some are already well known, but
others less so. The elliptic curve $E/K$ and integer $n \ge 2$ 
remain fixed throughout.

\smallskip

\begin{tabular}{llll}
& & \!\!\!\!\! base object & \!\!\!\!\!  twisted object \\
1. & Torsor divisor class pairs & $(E,[n.0])$ & $(C,[D])$ \\
2. & $n$-coverings & $(E,[n])$ & $(C,\pi)$ \\
3. & Brauer-Severi diagrams 
& $ [E \to \PP^{n-1}]$ & $[C \to S]$ \\
4. & $E[n]$-torsors & $(E[n],+)$ & $(\Phi,\mu)$ \\
5. & Comm. extns. of $E[n]$ by $\Gm$ 
&  $\Gm \times E[n]$ & $\Lambda$  \\
6. & Theta groups & $\Theta_E$ & $\Theta$. 
\end{tabular}

\medskip

Each of these interpretations appears frequently in our 
work (except for the fourth, which is required only for the flex
algebra method). The first three interpretations depend on the 
elliptic curve $E$ in an essential way: indeed 
$E$ may be recovered as the Jacobian of $C$. The fourth and fifth
interpretations depend only on $E[n]$ as a Galois module
equipped with the Weil paring.
If $n$ is odd then, by Lemma~\ref{thetaodd}, the same is true of the 
sixth interpretation.
We now go through each of the six interpretations
in turn.  

\subsection{First interpretation: Torsor divisor class pairs}
\label{tdcp}

\begin{Definition} (i) A torsor divisor class pair $(C,[D])$
is a torsor $C$ under $E$ together with a $K$-rational
divisor class $[D]$ on $C$ of degree $n$. The rationality means that $D$ is
linearly equivalent, but not necessarily equal, to all its 
Galois conjugates. \\
(ii) An isomorphism of torsor divisor class pairs $(C_1,[D_1]) \isom 
(C_2,[D_2])$ is an isomorphism of torsors $\phi:C_1 \isom C_2$ 
with $\phi^* D_2 \sim D_1$.
\end{Definition}

The trivial (or base) torsor divisor class pair is $(E,[n.0])$ 
where $0$ is the identity on $E$.
We recall that two divisors on an 
elliptic curve are linearly equivalent if and only if they have 
the same degree and the same sum.

\begin{Lemma}
\label{lem1.8}
Every torsor divisor class pair is a twist of $(E,[n.0])$.
\end{Lemma}
\begin{Proof} Let $(C,[D])$ be a torsor divisor class pair.
We choose an isomorphism of torsors $\phi : C \isom E$ 
defined over $\Kbar$. We then compose with a translation so 
that $\phi^* (n.0) \sim D$. 
\end{Proof}

\begin{Lemma} 
\label{auttdcp}
$\Aut(E,[n.0]) \isom E[n]$.
\end{Lemma}
\begin{Proof}
The automorphisms of $E$ as a torsor under itself are the translation
maps $\tau_P$ for $P \in E$. It suffices to note that 
$\tau_P^*(n.0) \sim n.0$ if and only if $P \in E[n]$.
\end{Proof}

By the twisting principle we obtain

\begin{Proposition} The torsor divisor class pairs, 
viewed as twists of $(E,[n.0])$, 
are parametrised up to isomorphism by $H^1(K,E[n])$. 
\end{Proposition}

\begin{Remark}
\label{kumCD}
The maps in the Kummer exact sequence
$$ E(K) \stackrel{\delta}{\ra} H^1(K,E[n]) \stackrel{\iota}{\ra} 
H^1(K,E) $$
are given by $\delta(P)= (E,[(n-1).0+P])$  and $\iota(C,[D]) = C$.
\end{Remark}

\subsection{Second Interpretation: $n$-coverings}

\begin{Definition} (i) A covering of $E$ is a pair $(C,\pi)$
where $C$ is a smooth projective curve and $\pi: C \to E$ 
is a non-constant morphism. \\
(ii) An isomorphism of coverings $(C_1,\pi_1) \isom (C_2,\pi_2)$ 
is an isomorphism of curves $\phi:C_1 \isom C_2$ with $\pi_1= \pi_2
\circ \phi$. 
\end{Definition}

We write $[n]$ for the multiplication-by-$n$ map on $E$.
The trivial (or base) $n$-covering of $E$ is $(E,[n])$.

\begin{Definition} An $n$-covering $(C,\pi)$ is a twist of $(E,[n])$.
\end{Definition}

\begin{Lemma} $\Aut(E,[n]) \isom E[n]$.
\end{Lemma}

\begin{Proof} Let $\phi:E \to E$ be an automorphism of $(E,[n])$.
Then we have $[n] = [n] \circ \phi$ and so $[n]\circ(\phi - 1)=0$.
We deduce that $\phi - 1$ is not surjective, and therefore constant.
It follows that $\phi$ is 
translation by an $n$-torsion point.
\end{Proof}

By the twisting principle we obtain

\begin{Proposition} The $n$-coverings of $E$ are 
parametrised up to isomorphism by $H^1(K,E[n])$. 
\end{Proposition}

\begin{Remark}
Given $(C,[D])$ a torsor divisor class pair, the corresponding
$n$-covering is $(C,\pi)$ where $\pi: C \to \Pic^0(C) \isom E$ 
is the map $P \mapsto [n.P-D]$. Conversely, given an $n$-covering
$(C,\pi)$ there exists an isomorphism $\phi:C \to E$
defined over $\Kbar$ making the diagram 
 $$ \xymatrix{ C \ar[dr]^{\pi} \ar[d]_\phi \\ E \ar[r]_{[n]} & E } $$
commute. We give $C$ the structure of torsor under $E$ via
$$(P,Q) \mapsto \phi^{-1}(P+\phi(Q)).$$ 
This definition is independent of the choice of $\phi$. 
The corresponding torsor divisor class pair is $(C, [\phi^*(n.0)])$.
The maps in the Kummer exact sequence
$$ E(K) \stackrel{\delta}{\ra} H^1(K,E[n]) \stackrel{\iota}{\ra} 
H^1(K,E) $$
are given by $\delta(P)= (E,\tau_P \circ [n])$  and $\iota(C,\pi) = C$.
\end{Remark}

\subsection{Third interpretation: Brauer-Severi diagrams}
\label{BSdiag}

\begin{Definition}
(i) A diagram $[C \to S]$ is a morphism from a torsor $C$ under $E$ 
to a variety $S$. \\
(ii) An isomorphism of diagrams $[C_1 \to S_1] \isom [C_2 \to S_2]$
is an isomorphism of torsors $\phi:C_1 \isom C_2$ together with
an isomorphism of varieties $\psi:S_1 \isom S_2$ making
the diagram
\begin{equation*}
 \xymatrix{
C_1 \ar[r] \ar[d]_{\phi} & S_1 \ar[d]^{\psi}   \\
C_2 \ar[r] & S_2  } 
\end{equation*}
commute.
\end{Definition}

The trivial (or base) diagram $[E \to \PP^{n-1}]$ is that
determined by the complete linear system $|n.0|$. 
We recall that a twist of projective space is 
called a Brauer-Severi variety.

\begin{Definition}
A Brauer-Severi diagram $[C \to S]$ is a twist of 
$[E \to \PP^{n-1}]$. 
In particular $S$ is a Brauer-Severi variety.
\end{Definition}

\begin{Lemma} $\Aut [E \to \PP^{n-1}] \isom E[n]$.
\label{autBSdiag}
\end{Lemma}
\begin{Proof} An automorphism $\phi$ of $E$ 
extends to an automorphism 
of $\PP^{n-1}$ if and only if $\phi^*(n.0) \sim n.0$. 
We are done by Lemma~\ref{auttdcp}.
\end{Proof}

By the twisting principle we obtain

\begin{Proposition} The Brauer-Severi diagrams
are parametrised up to isomorphism by $H^1(K,E[n])$. 
\end{Proposition}

\begin{Remark} Given a torsor divisor class pair $(C,[D])$
the complete linear system $|D|$ may be identified with the dual of a
Brauer-Severi variety $S$. There is then a natural morphism $C \to S$.
Conversely, given a Brauer-Severi diagram $[C \to S]$ we have 
$S \isom \PP^{n-1}$ over $\Kbar$. We pull back the hyperplane section 
on $\PP^{n-1}$ to give a $K$-rational divisor class $[D]$ on $C$. 
\end{Remark}

If $n \ge 3$ then the morphism $C \to \PP^{n-1}$ determined by
a complete linear system of degree $n$ is an embedding. The
image is called a genus one normal curve of degree $n$. In the case $n=3$
we get a smooth plane cubic. The homogeneous ideal of a genus one 
normal curve of degree $n \ge 4$ is generated by a vector space 
of quadrics of dimension $n(n-3)/2$.
The term ``elliptic normal curve'' is standard in the geometric 
literature, the word ``normal'' referring to the fact that the homogeneous
co-ordinate ring is integrally closed. We say ``genus one normal curve''
since we do not wish to imply that our curves have rational points.

\subsection{Fourth interpretation: $E[n]$-torsors}
\label{Entorsors}

\begin{Definition}
(i) An $E[n]$-torsor is a pair $(\Phi,\mu)$ where $\Phi$ 
is a zero-dimensional variety and $\mu:E[n] \times \Phi \to \Phi$ 
is a morphism that induces a simple transitive action on $\Kbar$-points. \\
(ii) An isomorphism of $E[n]$-torsors $(\Phi_1,\mu_1) \isom (\Phi_2,\mu_2)$ 
is an isomorphism of varieties $\Phi_1 \isom \Phi_2$
that respects the action of $E[n]$.
\end{Definition}

The trivial $E[n]$-torsor is $(E[n],+)$ where $+$ is the restriction
of the group law on $E$. In a manner entirely analogous to the proof
of Proposition~\ref{prototype} we obtain 

\begin{Proposition} The $E[n]$-torsors, viewed as twists of $(E[n],+)$,
are parametrised up to isomorphism by $H^1(K,E[n])$. 
\end{Proposition}

\begin{Remark} 
\label{remflex}
(i) Given a torsor divisor class pair $(C,[D])$ the
corresponding $E[n]$-torsor is the set of ``flex points'', i.e. 
$$\{ P \in C : n.P \sim D \}.$$ 
(ii) An $n$-covering $\pi: C \to E$ determines an 
$E[n]$-torsor $\pi^{-1}(0)$. \\
(iii) The connecting map $E(K) \to H^1(K,E[n])$ sends 
$P \mapsto [n]^{-1}(P)$. 
\end{Remark}

\begin{Remark} Generically, the splitting field of $\Phi$ has Galois group
the affine general linear group, $\AGL(2,n)$, which sits in an exact sequence
$$ 0 \to (\Z/n\Z)^2 \to \AGL(2,n) \to \GL(2,n) \to 0. $$
In the case $n=2$ this reduces to the exact sequence 
$0 \to V_4 \to S_4 \to S_3 \to 0$, as studied in \cite{CrInv}.
\end{Remark}

\subsection{Fifth Interpretation: Commutative extensions of $E[n]$ by $\Gm$.}
\label{commextns}

The following definition is common to our fifth and sixth interpretations
of $H^1(K,E[n])$. 

\begin{Definition}  
\label{defcentext}
(i) A central extension of $E[n]$ by $\Gm$ is an exact sequence 
of group varieties
$$ 0 \ra \Gm \stackrel{\alpha}{\ra} \Lambda \stackrel{\beta}{\ra} E[n] 
\ra 0 $$
with $\Gm$ contained in the centre of $\Lambda$. \\
(ii) An isomorphism of central extensions $\Lambda_1 \isom \Lambda_2$ 
is an isomorphism of group varieties $\phi: \Lambda_1 \isom \Lambda_2$ 
making the diagram
\begin{equation*}
 \xymatrix{
0 \ar[r] & \Gm \ar[r]^{\alpha_1} \ar@{=}[d] & \Lambda_1 \ar[r]^{\beta_1} 
\ar[d]_{\phi} & E[n] \ar[r] \ar@{=}[d] & 0  \\
0 \ar[r] & \Gm \ar[r]^{\alpha_2}& \Lambda_2  \ar[r]^{\beta_2} & E[n]
\ar[r] & 0  } 
\end{equation*}
commute.
\end{Definition}

We usually refer to $\Lambda$ as a central extension, the maps
$\alpha$ and $\beta$ being taken for granted. The trivial extension
is $\Lambda_0= \Gm \times E[n]$.

\begin{Lemma}
\label{extntwists} 
Every commutative extension of $E[n]$ by $\Gm$ 
is a twist of $\Lambda_0$. 
\end{Lemma}
\begin{proof}
Since $\Kbar^\times$ is a divisible group every commutative 
extension of $E[n]$ by $\Gm$ splits over $\Kbar$.
\end{proof}

\begin{Lemma}
\label{extnauts}
Let $\Lambda$ be any central extension of $E[n]$ by $\Gm$.
Then $\Aut(\Lambda) \isom \Hom(E[n],\Gm) \isom E[n]$.
\end{Lemma}
\begin{proof}
The automorphisms of $\Lambda$ take the form 
$x \mapsto \alpha(\pi (\beta (x))) x$
for $\pi: E[n] \to \Gm$ a homomorphism. This gives the first 
isomorphism. The second isomorphism comes from the Weil
pairing.
\end{proof}

By the twisting principle we obtain

\begin{Proposition} 
\label{extprop}
The commutative extensions of $E[n]$ by $\Gm$,
viewed as twists of $\Lambda_0$, are parametrised up to 
isomorphism by $H^1(K,E[n])$. 
\end{Proposition}

This result may be interpreted as giving an isomorphism
$$H^1(K,E[n]) \isom \Ext^1_K(E[n],\Gm)$$ 
cf. \cite{Milne}, Example 0.8.

\subsection{Sixth Interpretation: Theta groups}
\label{thetagps}

\begin{Definition}  
\label{deftheta}
(i) A theta group is a central extension of $E[n]$ by $\Gm$ 
$$ 0 \ra \Gm \stackrel{\alpha}{\ra} \Theta \stackrel{\beta}{\ra} E[n] \ra 0 $$
with commutator given by the Weil pairing, i.e.
$$ x y x^{-1} y^{-1} = \alpha (e_n(\beta x, \beta y)) $$
for all $x,y \in \Theta$.  \\
(ii) An isomorphism of theta groups is an isomorphism of central extensions
(see Definition~\ref{defcentext}).
\end{Definition}

In \S\ref{BSdiag} we considered the morphism 
$E \to \PP^{n-1}$ determined by the complete linear system $|n.0|$.
The action of $E[n]$ on $E$ by translation extends to an action on 
$\PP^{n-1}$ and so determines a map $\chi_E: E[n] \to \PGL_n$. 
Writing $\Theta_E$ for the inverse image of $\chi_E(E[n])$ in 
$\GL_n$ we obtain a commutative diagram with exact rows:
\begin{equation}
\label{ThetaEdiag}
 \xymatrix{
0 \ar[r] & \Gm \ar[r] \ar@{=}[d] & \Theta_E 
\ar[r] \ar[d] & E[n] 
\ar[r] \ar[d]^{\chi_E} & 0  \\
0 \ar[r] & \Gm \ar[r]      & \GL_n  \ar[r]    & \PGL_n
\ar[r]           & 0  } 
\end{equation}
It is clear that $\Theta_E$ is a central extension of
$E[n]$ by $\Gm$. To show it is a theta group, we must show
that it has commutator given by the Weil pairing.
In fact this may be used as the definition of the Weil pairing
(cf. Lemma~\ref{standardtheta}). 
For the relationship with the definition in 
Silverman \cite{Silverman}, Chapter III, \S8,
we refer to Mumford \cite{mumford}, \S\S20,23. 
There are issues of choice of sign here which we will ignore.

\begin{Lemma}
\label{thetaKbar} 
Every theta group is a twist of $\Theta_E$.
\end{Lemma}
\begin{Proof}
More generally we show that 
any two central extensions of $E[n]$ by $\Gm$, with the same commutator 
pairing, must necessarily be isomorphic over $\Kbar$.
Let $\Lambda_1$ and $\Lambda_2$ be two such extensions, and pick
a basis $S$, $T$ for $E[n]$.
Since $\Kbar^\times$ 
is a divisible group we may lift $S$, $T$ to elements 
$s_1,t_1 \in \Lambda_1$ and $s_2,t_2 \in \Lambda_2$ 
each of order $n$. Then there is an isomorphism 
$\Lambda_1 \isom \Lambda_2$, defined over $\Kbar$, 
uniquely determined by $s_1 \mapsto s_2$ and $t_1 \mapsto t_2$.
\end{Proof}

By Lemmas~\ref{extnauts},~\ref{thetaKbar} and the twisting
principle we obtain

\begin{Proposition} The theta groups for $E[n]$, viewed as twists 
of $\Theta_E$, are parametrised up to isomorphism by $H^1(K,E[n])$. 
\end{Proposition}

Let $(C,[D])$ be a torsor divisor class pair. We assume for simplicity
that $D$ is a $K$-rational divisor. Then the complete linear system
$|D|$ determines a morphism $C \to \PP^{n-1}$. 
The action of $E[n]$ on $C$ extends to an action on $\PP^{n-1}$
and so determines a map $\chi_C:E[n] \to \PGL_n$. We obtain a diagram
analogous to~(\ref{ThetaEdiag}):
\begin{equation}
\label{Thetadiag}
 \xymatrix{
0 \ar[r] & \Gm \ar[r] \ar@{=}[d] & \Theta 
\ar[r] \ar[d] & E[n] 
\ar[r] \ar[d]^{\chi_C} & 0  \\
0 \ar[r] & \Gm \ar[r]      & \GL_n  \ar[r]    & \PGL_n
\ar[r]           & 0  } 
\end{equation}
We show that this construction of $\Theta$ from $(C,[D])$ 
is compatible with our first and sixth interpretations of $H^1(K,E[n])$.
To do this, let $(C,[D])$ be the twist of $(E,[n.0])$ by 
$\xi \in H^1(K,E[n])$. Then there is an isomorphism of 
Brauer-Severi diagrams defined over $\Kbar$,
\begin{equation*}
 \xymatrix{
C \ar[r] \ar[d]_{\phi} & \PP^{n-1} \ar[d]^{\psi}   \\
E \ar[r]^{|n.0|} & \PP^{n-1}  } 
\end{equation*}
with $\sigma(\phi) \phi^{-1} = \tau_{\xi_\sigma}$. 
We lift $\psi$ to a matrix $B \in \GL_n$ so 
that conjugation by $B$ defines an isomorphism $\Psi: \Theta \isom \Theta_E$.
It is evident that $y_\sigma= \sigma(B) B^{-1}$ is an element of $\Theta_E$ 
projecting onto $\xi_\sigma$. Therefore
$$ \sigma(\Psi) \Psi^{-1} : \Theta_E \to \Theta_E \, ; \,\,
x \mapsto y_{\sigma} x y_{\sigma}^{-1} = e_n(\xi_\sigma,x) x $$
and $\Theta$ is the twist of $\Theta_E$ by $\xi$ as was to be shown.

In fact there is a more direct way to construct $\Theta$ from 
$(C,[D])$. We write $\tau_P : C \to C \, ; \,\, Q \mapsto \mu(P,Q)$ 
for the action of $P \in E$ on $C$.
\begin{Proposition}
\label{thetapairs}
Let $(C,[D])$ be a torsor divisor class pair with $D$ a $K$-rational 
divisor. Then the corresponding theta group is
$$ \Theta = \{ \, (f,T) \in \Kbar(C)^\times \times  E[n] \mid 
 \divv(f) = \tau_T^* D - D \, \}$$
with group law
\begin{equation}
\label{thetagrouplaw}
(f_1,T_1) * (f_2,T_2)= ( \tau_{T_2}^* (f_1) f_2 , T_1 + T_2 )
\end{equation}
and structure maps $\alpha: \lambda \mapsto (\lambda,0)$
and $\beta :(f,T) \mapsto T$.
\end{Proposition}
\begin{Proof}
The complete linear system $|D|$ determines a morphism
$$ C \to \PP( \CL(D)^*) \isom \PP^{n-1}$$ 
where $\CL(D)$ is the Riemann-Roch space
$$ \CL(D) = \{ f \in \Kbar(C)^\times \,\, |  
\, \divv (f)+D \ge 0 \} \cup \{ 0 \}. $$
Let $\iota: \Theta \to \End(\CL(D)^*)$ be given by
$$ \iota(f,T)(x) = (h \mapsto x(f \tau_T^*h)) $$
for all $x \in \CL(D)^*$ and $h \in \CL(D)$.
It may be verified that the diagram
\begin{equation*}
 \xymatrix{
C \ar[r] \ar[d]_{\tau_T} &  \PP( \CL(D)^*) \ar[d]^{\iota(f,T)}   \\
C \ar[r]  & \PP( \CL(D)^*)  } 
\end{equation*}
commutes. 
The group law on $\Theta$ is then that required to make
$\iota$ a group homomorphism. 
\end{Proof}

\section{The obstruction map}
\label{sect2}

We continue to take $E$ an elliptic curve over $K$ and $n \ge 2$ an 
integer. We recall from \S\ref{tdcp} that $H^{1}(K,E[n])$ 
parametrises the torsor divisor class pairs $(C,[D])$.
Here $[D]$ is a $K$-rational divisor class on $C$.

For any smooth projective variety $X$ over $K$ there is an 
exact sequence of Galois modules 
$$ 0 \to \Kbar^\times \to \Kbar(X)^\times \to \Div(X) \to \Pic(X) \to 0. $$
Splitting into short exact sequences and taking Galois cohomology 
we obtain an exact sequence (see \cite{Lichtenbaum} for details)
\begin{equation}
\label{grothbrauer}
0 \to K^\times \to K(X)^\times \to \Div(X)^{G_K} \to \Pic(X)^{G_K} 
\stackrel{\delta_X}{\to} \Br(K). 
\end{equation}

We define the obstruction map
$$ \Ob: H^1(K,E[n]) \to \Br(K) \,; \,\, (C,[D]) \mapsto \delta_C ([D]). $$
From~(\ref{grothbrauer}) we obtain the fundamental
property of the obstruction map, namely that $D$ is linearly 
equivalent to a $K$-rational divisor if and only if $\Ob(C,[D])=0$. 
We also have 

\begin{Lemma} 
\label{soluble => trivial obstruction}
Let $(C,[D])$ be a torsor divisor class pair.
If $C(K) \not= \emptyset$ 
(equivalently $C \isom E$ over $K$) then $\Ob(C,[D])=0$.
\end{Lemma}
\proof{}
Let $P \in C(K)$. 
By Riemann-Roch $D-(n-1)P$ is
linearly equivalent to a unique point $Q \in C$. 
The uniqueness statement proves that $Q$ is $K$-rational. 
So $D \sim (n-1)P+Q$ and the latter is $K$-rational. It follows
that $\Ob(C,[D])=0$. Alternatively the lemma follows from Remark~\ref{kumCD}.
\endproof

We give an alternative description of the obstruction map.
The base theta group $\Theta_E$ was defined in \S\ref{thetagps}.

\begin{Proposition}
\label{thetaobs}
The obstruction map $$ \Ob: H^1(K,E[n]) \to H^2(K,\Gm) $$
is obtained as the connecting map of (non-abelian) 
Galois cohomology for the exact sequence
$$ 0 \to \Gm \to \Theta_E \to E[n] \to 0. $$
\end{Proposition}
\proof{}
The proof is by means of a cocycle calculation. We start with
a torsor divisor class pair $(C,[D])$.
Since $[D]$ is a $K$-rational divisor class, there exist
rational functions $h_\sigma \in \Kbar(C)^\times$ with $\divv (h_\sigma) =
\sigma D - D$ for all $\sigma \in G_K$. 
By definition of the obstruction map, $\Ob(C,[D])$ 
is represented by the cocycle
\begin{equation}
\label{obcocycle}
(\sigma, \tau) \mapsto \frac{\sigma(h_\tau) h_\sigma}{h_{\sigma \tau}}. 
\end{equation}
Let $\xi \in H^1(K,E[n])$ describe $(C,[D])$ as a twist of $(E,[n.0])$.
This means there is an isomorphism of torsors
$\phi: C \isom E$, defined over $\Kbar$, with $\phi^* (n.0)=D$ and 
$\sigma(\phi) \phi^{-1} = \tau_{\xi_\sigma}$.
We recall that $\Theta_E$ is the theta group determined by
$(E,[n.0])$. By Proposition~\ref{thetapairs} we identify
$$ \Theta_E = \{ \, (f,T) \in \Kbar(E)^\times \times  E[n] 
\mid \divv(f) = \tau_T^*(n.0) - n.0 \, \}.$$
We lift $\xi_\sigma \in E[n]$ to a pair $(f_\sigma, \xi_\sigma)
\in \Theta_E$.
Then
$$ \begin{array}{rcl}
\divv (\phi^* f_\sigma) & = & \phi^* \tau_{\xi_\sigma}^* (n.0) - \phi^*
(n.0) \\
                    & = & (\sigma \phi)^* (n.0) - \phi^* (n.0) \\
                    & = & \sigma D - D. 
\end{array} $$
Taking $h_\sigma = \phi^* f_\sigma$ in~(\ref{obcocycle}) we obtain
$$ \begin{array}{rcl}
\Ob(\xi)(\sigma,\tau) & = &  \sigma (\phi^* f_\tau) \,\, \phi^*
f_\sigma \,\, ( \phi^* f_{\sigma\tau})^{-1}  \\
& = &  \tau_{\xi_\sigma}^* (\sigma f_\tau) \,\, f_\sigma  
\,\, f_{\sigma \tau}^{-1} \\
& = &  (\tau_{\xi_\sigma}^* (\sigma f_\tau) \,\, f_\sigma,
\xi_{\sigma \tau}) * (f_{\sigma \tau},\xi_{\sigma \tau})^{-1} \\
& = & \sigma(f_\tau,\xi_\tau) * (f_\sigma,\xi_\sigma) * (f_{\sigma
\tau},\xi_{\sigma \tau})^{-1} 
\end{array} $$
where $*$ is the group law~(\ref{thetagrouplaw}).
We recognise this final expression as the connecting map of Galois cohomology.
\endproof

\begin{Remark} We may identify $\ker(\Ob) = H^1(K, \Theta_E)$.
\end{Remark}

It is well known that taking Galois cohomology of the exact sequence 
$$ 0 \to \Gm \to \GL_n \to \PGL_n \to 0 $$
gives an injection
$$ \Delta: H^1(K,\PGL_n) \inj \Br(K)[n]. $$
This leads to our third interpretation of the obstruction map.

\begin{Lemma} 
\label{Ob=chi_*}
Let $[E \to \PP^{n-1}]$ be the 
base Brauer-Severi diagram
and let $\chi_E: E[n] \to \PGL_n$ describe the action of $E[n]$ on $E$.
Then the obstruction map is
$$ \chi_{E,*}: H^1(K,E[n]) \to H^1(K, \PGL_n). $$
\end{Lemma}
\begin{Proof} (Taken from \cite{cathy2}.)
Taking Galois cohomology in~(\ref{ThetaEdiag}) 
we obtain a commutative diagram
$$ \xymatrix{
H^1(K,E[n]) \ar[r] \ar[d]_{\chi_{E,*}} & H^2(K,\Gm) \ar@{=}[d] \\
H^1(K,\PGL_n) \ar[r]^\Delta & H^2(K,\Gm) } $$
The lemma follows from the description of the obstruction map
given in Proposition~\ref{thetaobs}.
\end{Proof}

Finally we interpret the obstruction map as a forgetful map.

\begin{Corollary} 
\label{obforget}
The obstruction map $\Ob: H^1(K,E[n]) \to \Br(K)$ 
sends $[C \to S]$ to $S$. 
\end{Corollary}
\begin{Proof} 
Let $[C \to S]$ be the twist of $[E \to \PP^{n-1}]$ by 
$\xi \in H^1(K,E[n])$. Then there is an isomorphism of Brauer-Severi 
diagrams defined over~$\Kbar$
\begin{equation*}
 \xymatrix{
C \ar[r] \ar[d]_{\phi} & S \ar[d]^{\psi}   \\
E \ar[r]^{} & \PP^{n-1}  } 
\end{equation*}
with $\sigma(\phi) \phi^{-1} = \tau_{\xi_\sigma}$. 
It follows that $\sigma(\psi) \psi^{-1} = \chi_E(\xi_\sigma)$ 
and so $S$ is the twist of $\PP^{n-1}$ by 
$\Ob(\xi)=\chi_{E,*}(\xi) \in H^1(K,\PGL_n)$.
\end{Proof}

If $[C \to S]$ is a Brauer-Severi diagram with $C(K) \not= \emptyset$
then clearly $S(K) \not= \emptyset$. So Corollary~\ref{obforget} gives an 
alternative proof of Lemma~\ref{soluble => trivial obstruction}.

\begin{Remark}
In general the obstruction map is not a group homomorphism.
But, as shown in~\cite{cathy2}, it is quadratic in the sense that \\
\indent (i) $\Ob(a \xi) = a^2 \Ob(\xi)$ for $a$ an integer, and \\
\indent (ii) $(\xi,\eta) \mapsto \Ob(\xi+\eta)-\Ob(\xi) - \Ob(\eta)$
is bilinear.
\end{Remark}

\section{The \'etale algebra}
\label{etal}

Let $R$ be the affine co-ordinate algebra of $E[n]$. 
It consists of all Galois equivariant maps from $E[n]$ to $\Kbar$. 
In symbols 
\begin{equation*}
R = \Map_K(E[n], \Kbar). 
\end{equation*}
For example, any rational function on $E$ defined over $K$ and not having
poles in~$E[n]$ will give an element of~$R$.
Since $E[n]$ is an \'etale $K$-scheme, $R$ is an \'etale algebra: it
is isomorphic to a product of (finite) field extensions of~$K$, one
for each $G_K$-orbit in~$E[n]$. If $T$ is a point in one such orbit,
then the corresponding field extension is $K(T)$.
If $n$ is prime then typically we have $R = K \times L$ where
$L/K$ is a field extension of degree $n^2-1$.

We also work with the algebra
$$ \Rbar = R \otimes_K \Kbar = \Map(E[n],\Kbar). $$
Note that the action of $\sigma \in G_K$ is 
$\alpha \mapsto (\sigma( \alpha):T \mapsto \sigma (\alpha (\sigma^{-1}T))).$
As a $\Kbar$-vector space $\Rbar$ has basis the $\delta_T$ 
for $T \in E[n]$ where
$$ \delta_{S}( T ) = \left\{ \begin{array}{ll} 1 & \text{ if } S=T, \\
0 & \text{ otherwise.} \end{array} \right. $$ 

In general, if $R_A$ and $R_B$ are the coordinate rings of two affine
$K$-schemes $A$ and~$B$, then $R_A \tensor_K R_B$ is the coordinate
ring of $A \times B$. 
So $R \tensor_K R$ is the algebra of Galois 
equivariant maps from $E[n] \times E[n]$ into~$\Kbar$, and 
$\Rbar \tensor_{\Kbar} \Rbar = (R \tensor_K R) \tensor_K \Kbar$
is the algebra of all such maps. 

The Weil pairing $ e_n: E[n] \times E[n] \to \mu_n$ determines an injection 
$$ w : E[n] \inj \Rbar^\times = \Map(E[n],\Kbar^\times) $$
via $w(S)(T)=e_n(S,T)$. 
We observe that the $w(T)$, for $T \in E[n]$, are not only maps, but
also homomorphisms $E[n] \to \Kbar^\times$. By the non-degeneracy of
the Weil pairing, all such homomorphisms arise in this way.
So if we define
$\partial: \Rbar^\times \to 
(\Rbar \tensor_{\Kbar} \Rbar)^\times $ via
\begin{equation}
\label{defpartial}
 (\partial \alpha)(T_1,T_2) = 
\frac{ \alpha(T_1) \alpha(T_2) }{\alpha(T_1+T_2)},
\end{equation}
then there is an exact sequence
\begin{equation}
\label{exseq1}
 0 \ra E[n] \stackrel{w}{\ra} \Rbar^\times 
          \stackrel{\partial}{\ra} (\Rbar \tensor \Rbar)^\times.
\end{equation}

By a generalised version of Hilbert's theorem 90 
(which reduces by Shapiro's lemma to the usual version 
of Hilbert's theorem 90 applied
to each constituent field of $R$) we have
$$ H^1(K,\Rbar^\times) = 0.$$

We use these observations to define group homomorphisms
$$ w_1 : H^1(K,E[n]) \to R^\times /(R^\times)^n $$
and 
$$ w_2 : H^1(K,E[n]) \to (R \otimes R)^\times/ \partial R^\times. $$
It is convenient to give both definitions at once. We start
with $\xi \in H^1(K,E[n])$ and use Hilbert's theorem 90 to write
$w(\xi_\sigma) = \sigma(\gamma)/\gamma$ for some $\gamma \in \Rbar^\times$.
Then $\alpha=\gamma^n$ and $\rho = \partial \gamma$ are Galois
invariant and so belong to $R^\times$ and $(R \otimes R)^\times$
respectively. We define 
$w_1(\xi)=\alpha (R^\times)^n$ and $w_2(\xi)= \rho \, \partial R^\times$. 
If we change $\xi$ by a coboundary, say $\sigma(T)-T$, then $\gamma$
is multiplied by $w(T)$. Since $w(T)^n=1$ and $\partial (w(T))=1$ this
leaves the values of $\alpha$ and $\rho$ unchanged. The only 
remaining freedom is to multiply $\gamma$ by an element of $R^\times$.
This has the effect of multiplying $\alpha$ and $\rho$ by 
elements of $(R^\times)^n$ and $\partial R^\times$ respectively.
It follows that $w_1$ and $w_2$ are well defined.

The map $w_1$ is in fact the composite 
$$ H^1(K,E[n]) \stackrel{w_*}{\ra} H^1(K,\mu_n(\Rbar)) 
\stackrel{k}{\ra} R^\times/(R^\times)^n $$
where $w_*$ is induced by $w$, and $k$ is the Kummer
isomorphism. 
\begin{Lemma} If $n$ is prime then $w_1$ is injective.
\end{Lemma}
\begin{Proof}
See \cite{DSS}, Proposition 7, or 
\cite{SchaeferStoll}, Corollary 5.1. 
\end{Proof}
In general $w_1$ is not injective. For example, taking
$n=4$ and $E/\Q$ the elliptic curve $y^2=x^3+x+2/13$, 
it may be shown that $w_1$ has kernel of order 2.

\begin{Lemma}
\label{w2inj}
The map $w_2$ is injective.
\end{Lemma}
\begin{Proof}
Let $\xi$ belong to the kernel of $w_2$. Then $w(\xi_\sigma)=
\sigma(\gamma)/\gamma$ for some $\gamma \in \Rbar^\times$.
Multiplying $\gamma$ by an element of $R^\times$ we may suppose that
$\partial \gamma =1$. Then~(\ref{exseq1}) gives $\gamma=w(T)$
for some $T \in E[n]$. Since $w$ is injective it follows that
$\xi_\sigma= \sigma(T)-T$. Hence $\xi$ is a coboundary.
\end{Proof}

In \S\ref{commextns} we showed that $H^1(K,E[n])$ parametrises
the commutative extensions of $E[n]$ by $\Gm$. This point of
view will help us determine the image of $w_2$.
By Hilbert's theorem 90 every central extension
$$ 0 \to \Gm \to \Lambda \to E[n] \to 0 $$
has a Galois equivariant section $\phi : E[n] \to \Lambda$.
In general $\phi$ is not a group homomorphism. 
The possible choices of $\phi$ differ by elements of 
$$ \Map_K(E[n],\Kbar^\times) = R^\times. $$

We define the first and second invariants of $\Lambda$. The first
invariant is $\inv_1(\Lambda) = \alpha (R^\times)^n$ where 
$\alpha \in R^\times$ satisfies
\begin{equation}
\label{alphadef}
 \phi(T)^n = \alpha(T)
\end{equation}
for all $T \in E[n]$. The second invariant is $\inv_2(\Lambda) = \rho 
\, \partial R^\times$  where $\rho \in (R \otimes R)^\times$ satisfies
\begin{equation}
\label{rhodef}
 \phi(T_1) \phi(T_2)  = \rho(T_1,T_2) \phi(T_1+T_2). 
\end{equation}
for all $T_1,T_2 \in E[n]$.
Notice that $\inv_1(\Lambda)$ and $\inv_2(\Lambda)$ 
depend only on $\Lambda$ and not on the choice of section $\phi$.

\begin{Lemma} 
\label{relatetocommext}
Let $\Lambda$ be the twist of $\Lambda_0$ 
by $\xi \in H^1(K,E[n])$. Then
$ \inv_1(\Lambda) = w_1(\xi)$ and  $\inv_2(\Lambda) = w_2(\xi)$.
\end{Lemma}

\begin{Proof}
By hypothesis there is a commutative diagram
$$
 \xymatrix{
0 \ar[r] & \Gm \ar[r] \ar@{=}[d] & \Lambda \ar[r] 
\ar[d]_{\psi} & E[n] \ar[r] \ar@{=}[d] & 0  \\
0 \ar[r] & \Gm \ar[r] & \Lambda_0  \ar[r] & E[n]
\ar[r] & 0  } 
$$
with $\sigma(\psi) \psi^{-1} : x \mapsto e_n(\xi_\sigma,x)x$.
We use Hilbert's theorem 90 to write $w(\xi_\sigma)=\sigma(\gamma)/\gamma$
for some $\gamma \in \Rbar^\times$. Let $\phi_0: E[n] \to \Lambda_0$
be the natural section for $\Lambda_0$. 
A calculation reveals that
$$ \phi : E[n] \to \Lambda \,; \,\, 
T \mapsto \gamma(T) \psi^{-1}(\phi_0(T)) $$
is a Galois equivariant section for $\Lambda$. So by~(\ref{alphadef}) 
and~(\ref{rhodef}) we have $\alpha=\gamma^n$ and $\rho= \partial \gamma$
as required. 
\end{Proof}

\begin{Remark}
It is clear that a central extension of $E[n]$ by $\Gm$ is uniquely
determined up to isomorphism by its second invariant. Thus 
Lemma~\ref{relatetocommext} gives an alternative 
proof of Lemma~\ref{w2inj}.
\end{Remark}

We extend~(\ref{exseq1}) to a complex
  \[ 0 \ra E[n] \stackrel{w}{\ra} \Rbar^\times
       \stackrel{\partial}{\ra} (\Rbar \tensor \Rbar)^\times
       \stackrel{\partial}{\ra} 
        (\Rbar \tensor \Rbar \tensor \Rbar)^\times
  \]
where the second $\partial$ is given by  
$$ (\partial \rho)(T_1,T_2,T_3) = \frac{ \rho(T_1,T_2) \rho(T_1+T_2,T_3) }
{\rho(T_1,T_2+T_3) \rho(T_2,T_3) }. $$ 
For each  $\rho \in (R \tensor R)^\times$ we write 
$\rho^\op$ for the element obtained by switching the operands, {\em i.e.}
$\rho^\op(T_1,T_2) = \rho(T_2,T_1)$.

\begin{Lemma}
\label{im=H}
The image of $w_2$ is $$H= \{ \, \rho \in (R \otimes R)^\times \mid 
\rho = \rho^\op \text{ and }  \partial \rho = 1 \} / \partial R^\times.$$
\end{Lemma}

\begin{Proof} The conditions  $\rho = \rho^\op$  and $\partial \rho = 1$  
express the fact that $\Lambda$ is commutative and associative. 
Conversely, if $\rho \, \partial R^\times \in H$ then we define a new
multiplication on $\Gm \times E[n]$ via
$$ (\lambda_1, T_1) * (\lambda_2,T_2)  = (\lambda_1 \lambda_2
\rho(T_1,T_2) ,T_1+T_2). $$
This gives the required commutative extension of $E[n]$ by $\Gm$. 
\end{Proof}

\begin{Corollary}
\label{kbarcor}
If $\rho \, \partial R^\times \in H$ then $\rho = \partial \gamma$ for some
$\gamma \in \Rbar^\times$.
\end{Corollary}
\begin{Proof}
This is the case $K=\Kbar$ of the last lemma.
\end{Proof}

\begin{Remark}
Lemma~\ref{im=H} may equally be deduced from Corollary~\ref{kbarcor}
by taking Galois cohomology of the short exact sequence 
$$ 0 \ra E[n] \stackrel{w}{\ra} \Rbar^\times \stackrel{\partial}{\ra}
\ker(\partial | \Sym^2(\Rbar)^\times ) \ra 0 $$
where $\Sym^2(\Rbar) = \{ \rho \in \Rbar \otimes 
\Rbar \mid \rho = \rho^\op \}$.
\end{Remark}

In applications we take $K$ a number field and $n$ a prime.
As explained in the introduction, we compute $\Sel^{(n)}(E/K)$
first as a subgroup of $R^\times/(R^\times)^n$ and then convert it to a 
subgroup of $H \subset (R \otimes R)^\times/\partial R^\times$. To make
this conversion, we let $\kappa$ be the map making the diagram
$$ \xymatrix{ H^1(K,E[n]) \ar[d]_-{w_2} \ar[dr]^-{w_1} & \\
H \ar[r]^-{\kappa}   & R^\times/(R^\times)^n } $$
commute. If $\rho \partial R^\times \in H$ then by Corollary~\ref{kbarcor} 
we have $\rho = \partial \gamma$ for some $\gamma \in \Rbar^\times$. 
Tracing through the definitions we find that $\kappa(\rho \partial R^\times)
= \alpha (R^\times)^n$ where $\alpha = \gamma^n \in R^\times$. 

\begin{Lemma}
Let $\alpha (R^\times)^n$ belong to the image of $w_1$. Then there
exists $\rho \in \Sym^2(R)^\times$ with
(i) $\partial \alpha = \rho^n$,
(ii) $\alpha(T) = \prod_{i=0}^{n-1} \rho(T,iT)$ for all $T \in E[n]$,
and (iii) $\partial \rho = 1.$ 
Moreover if $\rho \in \Sym^2(R)^\times$ satisfies (ii) and (iii)
then $\kappa(\rho \, \partial R^\times) = \alpha (R^\times)^n$.
\end{Lemma}
\begin{Proof}
Since $w_1= \kappa \circ w_2$ there exists $\rho \, \partial R^\times \in H$
with $\kappa(\rho \, \partial R^\times) = \alpha (R^\times)^n$. 
Corollary~\ref{kbarcor} gives $\rho = \partial \gamma$ 
for some $\gamma \in \Rbar^\times$. Multiplying
$\gamma$ by a suitable element of $R^\times$ we may suppose 
that $\gamma^n=\alpha$. Conditions (i), (ii) and (iii) follow at once. 

Conversely if $\rho \in \Sym^2(R)^\times$ satisfies (ii) and (iii)
then $\rho \, \partial R^\times \in H$ and so $\rho = \partial \gamma$
for some $\gamma \in \Rbar^\times$. By (ii) we deduce $\alpha= \gamma^n$.
\end{Proof}

\begin{Remark}
If $\Sym^2(R)$ contains no non-trivial $n$th roots of unity,
then we construct $\rho$ from $\alpha$ by taking the unique
$n$th root of $\partial \alpha$. There is no need to check
conditions (ii) and (iii).
\end{Remark}

We have identified $H^1(K,E[n])$ with a subgroup 
$H \subset (R \otimes R)^\times/\partial R^\times$. 
We did this using commutative extensions of $E[n]$ by $\Gm$. 
But we may equally work with theta groups (cf. \S\ref{thetagps}).

\begin{Lemma} 
\label{thetarho}
Let $\Theta$ be the twist of $\Theta_E$
by $\xi \in H^1(K,E[n])$. Then
$\inv_1(\Theta) =  w_1(\xi) \inv_1(\Theta_E) $ 
and  $\inv_2(\Theta) = w_2(\xi) \inv_2(\Theta_E)$.
\end{Lemma}
\begin{Proof} 
The proof is similar to that of Lemma~\ref{relatetocommext}.
\end{Proof}

Let $\eps \in (R \otimes R)^\times$ with 
$\inv_2(\Theta_E) = \eps \, \partial R^\times$. Then the 
theta groups for $E[n]$, or rather their second invariants, 
make up the coset
$$ \eps H= \{ \rho \in (R \otimes R)^\times \mid 
\rho (\rho^\op)^{-1} =e \text{ and } \partial \rho = 1 \} 
/ \partial R^\times$$
where $e \in \mu_n(R \otimes R)$ is the Weil pairing.
We discuss several ways of computing $\eps$.
Our first method uses the definition of $\Theta_E$ as a subgroup of
$\GL_n$. 
Recall that we mapped  $E \to \PP^{n-1}$ via $|n.0|$. 
Let $M$ be a matrix in 
$$ \GL_n(R) = \Map_K (E[n], \GL_n(\Kbar)) $$
that describes the action of $E[n]$ on $\PP^{n-1}$.
Then $\Theta_E$ has Galois equivariant section $T \mapsto M_T$ 
and so $\eps \in (R \otimes R)^\times$ is determined by 
$$ M_{T_1} M_{T_2} = \eps(T_1,T_2) M_{T_1+T_2}. $$ 

Our second method uses the description of $\Theta_E$ 
obtained by taking $(C,[D])=(E,[n.0])$ in Proposition~\ref{thetapairs}.
Let $F$ be a rational function in
$$ R(E)^\times = \Map_K(E[n], \Kbar(E)^\times) $$
with $\divv(F_T)=n.T-n.0$. Then $\Theta_E$ has Galois-equivariant 
section $T \mapsto (F_T,-T)^{-1}$.
Using the group law~(\ref{thetagrouplaw}) we obtain
$$ \eps(T_1,T_2) = \frac{ F_{T_1+T_2} (P)}{F_{T_1}(P)F_{T_2}(P-T_1)}$$
where the righthand side is constant as a function of $P \in E$.

\medskip

If $n$ is odd then we are spared the above calculations.
\begin{Lemma} 
\label{thetaodd}
If $n$ is odd, say $n=2m-1$, then $\inv_1(\Theta_E)$ is trivial 
and $\inv_2(\Theta_E)= e^m \partial R^\times$. In particular $\Theta_E$
depends on $E[n]$ and the Weil pairing, but not on $E$.
\end{Lemma}
\begin{Proof}
As before we map $E \to \PP^{n-1}$ via 
$|n.0|$. The action of $E[n]$ on $E$ determines 
$\chi_E:E[n] \to \PGL_n$. Likewise the negation map $[-1]$ 
on $E$ determines an element $\iota \in \PGL_n$. We claim that
for each torsion point $T \in E[n]$ there is a unique lift 
$M_T$ of $\chi(T)$ to $\GL_n$ such that (i) $\iota M_T \iota^{-1} = M_T^{-1}$ 
and (ii) $M_T^n=I$. Indeed the first condition determines $M_T$ up
to sign, and implies $M_T^n= \pm I$. Then, since $n$ is odd, 
the second condition determines a unique choice of this sign.

The uniqueness statement tells us that the map $\phi:T \mapsto M_T$
is Galois equivariant. A short calculation 
(using (i), (ii) and the commutator condition) reveals 
that $M_S M_T = e_n(S,T)^m M_{S+T}$ for all $S,T \in E[n]$.
Substituting in~(\ref{alphadef}) and~(\ref{rhodef}) we get
$\alpha=1$ and $\rho=e^m$ as required.
\end{Proof}

If $n$ is odd then the restriction of $w_1$ to the kernel of 
the obstruction map has the following alternative interpretation.

\begin{Corollary}
\label{detMcor}
Assume $n$ is odd. 
Let $[C \to \PP^{n-1}]$ be the Brauer-Severi diagram 
determined by $\xi \in H^1(K,E[n])$ and let $M \in \GL_n(R)$ describe the 
action of $E[n]$ on $C$. Then $w_1(\xi) = (\det M) (R^\times)^n$. 
\end{Corollary}
\begin{Proof} 
Let $\Theta$ be the theta group determined by $[C \to \PP^{n-1}]$.
Then $\Theta$ is the twist of $\Theta_E$ by $\xi \in H^1(K,E[n])$. 
By Lemmas~\ref{thetarho} and~\ref{thetaodd} we have 
$\inv_1(\Theta) = w_1(\xi)$. Therefore $w_1(\xi) = \alpha (R^\times)^n$
where $\alpha \in R^\times$ is determined by $M^n = \alpha I_n$.
The next lemma shows that if $T \in E[n]$ has order $r$ then 
$M_T$ has characteristic polynomial of the form $(X^r-c)^{n/r}$. 
We deduce that $\det(M) = \alpha$ as required.
\end{Proof} 

In the following lemma we assume that $K$ is algebraically
closed, and fix $\zeta_n \in K$ a primitive $n$th root of unity.
\begin{Lemma} 
\label{standardtheta}
Let $C \subset \PP^{n-1}$ be a genus one normal curve
with Jacobian $E$. Let $T_1$, $T_2$  be a basis for $E[n]$ with 
$e_n(T_1,T_2)= \zeta_n$. Then we can choose co-ordinates on $\PP^{n-1}$ 
so that $T_1$, $T_2$ act on $C$ via
$$ M_1 = \begin{pmatrix}
    1   &    0   &    0    & \cdots &      0     \\
    0   & \zeta_n  &    0    & \cdots &      0     \\
    0   &    0   & \zeta_n^2 & \cdots &      0     \\ 
 \vdots & \vdots & \vdots  &        &   \vdots   \\
    0   &   0    &    0    & \cdots & \zeta_n^{n-1}  
\end{pmatrix}, \quad  
    M_2 = \begin{pmatrix}
   0   &    0   & \cdots &   0    &      1     \\ 
    1   &    0   & \cdots &   0    &      0     \\
    0   &    1   & \cdots &   0    &      0     \\ 
 \vdots & \vdots &        & \vdots &   \vdots   \\
    0   &   0    & \cdots &   1    &      0       
\end{pmatrix}. $$
\end{Lemma}
\begin{Proof} 
This is quite standard. See for example \cite{cathy}.
\end{Proof}

\section{From extensions to enveloping algebras}
\label{algebras}

\subsection{Enveloping algebras}

We consider $K$-algebras $A$ that are finite dimensional over $K$.
The unit group 
of $\Abar = A \otimes_K \Kbar$ may be viewed as 
(the $\Kbar$-rational points of) a $K$-group variety.
For instance, if $A$ is the matrix algebra $\Mat_n(K)$ then this
construction yields $\GL_n$.

\begin{Definition}
\label{defemb}
Let $\Lambda$ be a central extension of $E[n]$ by $\Gm$.
Let $A$ be a $K$-algebra with $[A:K]=n^2$.
An {\em embedding} of $\Lambda$ in $A$ is a 
morphism of  $K$-group varieties
$\iota : \Lambda \to \Abar^\times$
such that \\
(i) $\iota$ preserves scalars, {\em i.e.} $\iota(\lambda)=\lambda$ 
for all $\lambda \in \Kbar^\times$, \\
(ii) the image of $\iota$ spans $\Abar$ as a $\Kbar$-vector space. 
\end{Definition}

\begin{Lemma}
\label{embexists}
Every central extension of $E[n]$ by $\Gm$ embeds in a $K$-algebra.
\end{Lemma}
\begin{Proof}
Let $\phi: E[n] \to \Lambda$ be a Galois equivariant 
section for $\Lambda$. As in \S\ref{etal}
we write $R$ for the \'etale algebra of $E[n]$ and recall that
$\Rbar$ is a $\Kbar$-vector space with basis the $\delta_T$ for $T \in E[n]$.
We define a Galois equivariant inclusion of $\Lambda$ in $\Rbar$ via
$$ \lambda \phi(T) \mapsto \lambda \delta_T $$
for all $\lambda \in \Kbar^\times$ and $T \in E[n]$. 
The group law on $\Lambda$ extends uniquely
to a new $\Kbar$-algebra multiplication on $\Rbar$, 
which in turn descends to a $K$-algebra multiplication
$$ * : R \times R \to R. $$
Thus $\Lambda$ embeds in the $K$-algebra $A = (R,+,*)$. 
\end{Proof}

\begin{Lemma}
\label{natural}
Let $\Lambda_1$ and $\Lambda_2$ be central extensions of $E[n]$ by $\Gm$ 
embedding in $K$-algebras $A_1$ and $A_2$. 
Then every isomorphism of central extensions 
$\psi: \Lambda_1 \isom \Lambda_2$ extends uniquely to an isomorphism 
of $K$-algebras $\Psi:A_1 \isom A_2$.
\end{Lemma}

\begin{Proof}
We construct $\Psi$ from $\psi$ by extending linearly to
an isomorphism of $\Kbar$-algebras, and then restricting to
$K$-algebras. 
Condition (ii) of Definition~\ref{defemb} ensures that $\Psi$
is unique.
\end{Proof}

\begin{Definition} 
Let $\Lambda$ be a central extension of $E[n]$ by $\Gm$.
If $\Lambda$ embeds in a $K$-algebra $A$ then 
$A$ is the {\em enveloping algebra} of $\Lambda$. 
\end{Definition}

We have shown that enveloping algebras exist and are unique up to 
isomorphism. Next we outline a method for computing them.
The addition law
$$ E[n] \times E[n] \to E[n] $$
gives rise to the comultiplication
$$ \Delta: R \to R \otimes R $$
with $\Delta(\alpha)(T_1,T_2)=\alpha(T_1+T_2)$.
Viewing $R \otimes R$ as an $R$-algebra via $\Delta$
there is a trace map 
$$ \Tr : R \otimes R \to R. $$
In terms of functions it is given by
\begin{equation}
\label{tracedef}
\Tr(\rho)(T) = \sum_{T_1+T_2=T} \rho(T_1,T_2). 
\end{equation}
It may also be built out of the trace maps for the
constituent fields of $R \otimes R$ and $R$.

\begin{Lemma} 
\label{mrho}
Let $\Lambda$ be a central extension of $E[n]$ by $\Gm$. If 
$\inv_2(\Lambda) = \rho \, \partial R^\times$ then $\Lambda$ has 
enveloping algebra $(R,+,*_\rho)$ where
$$  z_1 *_\rho z_2= \Tr( \rho . z_1 \otimes z_2). $$
\end{Lemma}
\begin{Proof} By hypothesis there exists $T \mapsto \phi(T)$ 
a Galois equivariant section for $\Lambda \to E[n]$ with
$$ \phi(T_1) \phi(T_2) =  \rho(T_1,T_2) \phi(T_1+T_2). $$ 
Following the proof of Lemma~\ref{embexists} we obtain 
$$ \delta_S *_\rho \delta_T =  \rho(S,T) \delta_{S+T}. $$
We write this multiplication in a form that descends to $K$: 
$$ \begin{array}{rcl}
z_1 *_\rho z_2 & = & (\sum_T z_{1}(T) \delta_{T}) *_\rho 
(\sum_T z_{2}(T) \delta_{T}) \\
& = & \sum_T ( \sum_{T_1+T_2=T} \rho(T_1,T_2) 
z_1(T_1) z_2(T_2)) \delta_{T} \\
& = & \Tr(\rho.z_1 \tensor z_2). 
\end{array} $$
\end{Proof}

Lemma~\ref{natural} tells us that if $\Lambda_1$ and $\Lambda_2$
are isomorphic (as central extensions) then $A_1$ and $A_2$ are 
isomorphic (as $K$-algebras). More concretely we have

\begin{Lemma}
\label{concrete}
Let $A_1 = (R,+, *_{\rho_1})$ and  $A_2 = (R,+, *_{\rho_2})$
be the enveloping algebras determined by 
$\rho_1, \rho_2 \in (R \otimes R)^\times$ with 
$\partial \rho_1= \partial \rho_2 = 1$. If $\rho_1= \rho_2 \partial \gamma$
for some $\gamma \in R^\times$ then there is an isomorphism of 
$K$-algebras 
$$ A_1 \isom  A_2 \, ; \,\, z \mapsto 
\gamma. z $$ 
where the multiplication is that in $R$.
\end{Lemma}

\begin{Proof} We compute
$$ \begin{array}{rcl}
\gamma. (z_1 *_{\rho_1} z_2) 
&=& (T \mapsto \gamma(T)\sum_{T_1+T_2=T} \rho_1(T_1,T_2) z_1(T_1)z_2(T_2)) \\
&=& (T \mapsto \sum_{T_1+T_2=T} \rho_2(T_1,T_2) 
\gamma(T_1)  z_1(T_1) \gamma(T_2) z_2(T_2)) \\
& = & (\gamma.z_1) *_{\rho_2} (\gamma.z_2)  
\end{array} $$
\end{Proof}

\subsection{The flex algebra}

Let $\Phi$ be an $E[n]$-torsor and let $F$ be the \'etale 
algebra of $\Phi$, {\em i.e.}
$$ F = \Map_K(\Phi,\Kbar). $$
Since $E[n]$ acts on $\Phi$ it also acts on 
$\Fbar^\times = \Map(\Phi,\Kbar^\times)$. 
The ``eigenvectors'' for this action form a group
$$ \Lambda = \left\{ z \in \Fbar^\times  \Bigg| \begin{array}{c} 
 \text{ there exists } T \in E[n] 
\text{ such that } \\ z(S+P)= e_n(S,T) z(P) \\ \text{ for all } S \in E[n]
\text{ and } P \in \Phi \end{array} \right\}. $$
Thus we obtain a commutative extension
$$ 0 \ra \Gm \stackrel{\alpha}{\ra} \Lambda \stackrel{\beta}{\ra} E[n] \ra 0 $$
where in the above notation $\beta(z)=T$. 
We show that this construction of $\Lambda$ from $F$ is compatible
with our fourth and fifth interpretations of $H^1(K,E[n])$. 
To do this let $\Phi$ be the twist of $E[n]$ by
$\xi \in H^1(K,E[n])$. This means there is an isomorphism of
$E[n]$-torsors $\psi: \Phi \to E[n]$, defined over $\Kbar$,
with $\sigma(\psi)\psi^{-1} = \tau_{\xi_\sigma}$. We define
$$ \psi^* : \Fbar= \Map(\Phi,\Kbar) \to \Map(E[n],\Kbar)=\Rbar $$
via $\psi^*(z)(P)=z(\psi^{-1}(P))$. Then $\psi^*$ restricts
to an isomorphism of central extensions $\gamma: \Lambda \isom \Lambda_0$ where
$$ \Lambda_0 = \Gm \times E[n] = \{ \, \lambda w(T) \in \Rbar^\times 
\mid \lambda \in \Kbar^\times, \,\, T \in E[n] \, \}. $$
We find  
$$\sigma(\gamma)\gamma^{-1}: \Lambda_0 \to \Lambda_0 \, ; \,\, x \mapsto
e_n(\xi_\sigma,x)x. $$
So $\Lambda$ is the twist of $\Lambda_0$ by $\xi$ as was to be shown.
We summarise the above discussion in

\begin{Proposition} 
\label{flexenv}
The enveloping algebra of a commutative extension
of $E[n]$ by $\Gm$ is the \'etale algebra of the corresponding $E[n]$-torsor.
\end{Proposition}

If $\xi \in H^1(K,E[n])$ with $w_2(\xi)= \rho \partial R^\times$
then by Lemmas~\ref{relatetocommext} and~\ref{mrho} 
the flex algebra is $(R,+,*_\rho)$.
Taking $\rho=1$ should give the \'etale algebra of $E[n]$. 
We recognise $(R,+,*_1)$ 
as the group algebra of $E[n]$. It is isomorphic to $R$ 
via the Fourier transform
$\alpha \mapsto \alphahat$ where
$$ \alphahat(S) = \frac{1}{n^2} \sum_T e_n(S,T) \alpha(T). $$

\subsection{The obstruction algebra}
\label{4.3}

\begin{Lemma} 
\label{extrivob}
The base theta group $\Theta_E$ embeds in $\Mat_n(K)$. 
In particular $\Theta_E$ spans $\Mat_n(\Kbar)$ as a $\Kbar$-vector space.
\end{Lemma}
\begin{Proof}
We recall that $\Theta_E$ was defined 
in \S\ref{thetagps} as a subgroup of $\GL_n$.
Let $T \mapsto M_T$ be a section for
$\Theta_E \to E[n]$. By Definition~\ref{deftheta} we have
$$ M_S M_T M_S^{-1} M_T^{-1} = e_n(S,T) $$
for all $S,T \in E[n]$. From the non-degeneracy of the Weil pairing
it follows that the $M_T$ are linearly independent over 
$\Kbar$. A dimension count shows that they span $\Mat_n(\Kbar)$. 
This proves the second statement of the lemma. The first now follows
by Definition~\ref{defemb}.
\end{Proof}
 
\begin{Lemma} 
\label{A=csa}
If $\Theta$ is a theta group for $E[n]$ with enveloping algebra $A$
then $A$ is a central simple $K$-algebra with $[A:K]=n^2$. 
In particular $\Theta$ spans $\Abar$ as a $\Kbar$-vector space.
\end{Lemma}
\begin{Proof} We saw in Lemma~\ref{thetaKbar} that $\Theta$ is a twist of 
$\Theta_E$. It follows by Lemma~\ref{natural} that $A$ is a twist
of $\Mat_n(K)$. 
\end{Proof}

We use Definition~\ref{defemb} and Lemma~\ref{A=csa} 
to build a commutative diagram with exact rows
\begin{equation}
\label{diagram**}
 \xymatrix{
0 \ar[r] & \Gm \ar[r] \ar@{=}[d] & \Theta \ar[r] \ar[d] & E[n]
\ar[r] \ar[d] & 0  \\
0 \ar[r] & \Gm \ar[r]  & \Abar^\times  \ar[r]  & \Aut(A)
\ar[r]           & 0.  } 
\end{equation}
The second row is exact by the Noether-Skolem theorem.
The righthand vertical arrow is the composite of 
$E[n]= \Aut(\Theta)$, cf. Lemma~\ref{extnauts}, and
$\Aut(\Theta) \to \Aut(A)$, cf. Lemma~\ref{natural}.
We recognise~(\ref{ThetaEdiag}) and~(\ref{Thetadiag}) 
as special cases of this new diagram. 

\begin{Proposition}
\label{newob}
If $\Theta$ is a theta group for $E[n]$ with enveloping algebra $A$
then the obstruction map sends the class of $\Theta$ in $H^1(K,E[n])$ to the
class of $A$ in $\Br(K)$. 
\end{Proposition}
\proof{}
Let $\Theta_E$ be the base theta group, {\em i.e.} the theta group 
associated to $(E,[n.0])$. 
Let $\chi_E:E[n] \to \PGL_n$ be the right-hand vertical arrow 
in~(\ref{ThetaEdiag}). According to Proposition~\ref{Ob=chi_*} the 
obstruction map is  $$\chi_{E,*}:H^1(K,E[n]) \to H^1(K,\PGL_n).$$

Now let $\Theta$ be the twist of $\Theta_E$ by $\xi \in H^1(K,E[n])$. 
This means there is an isomorphism of theta groups 
$\gamma:\Theta \isom \Theta_E$, defined over $\Kbar$, with
$$ \sigma(\gamma) \gamma^{-1}: x \mapsto e_n(\xi_\sigma,x)x $$
for all $x \in \Theta_E$. 
It follows from the commutator condition 
that $\sigma(\gamma) \gamma^{-1}$ is conjugation by 
any lift of $\xi_\sigma$ to $\Theta_E$.
Then Lemma~\ref{natural} tells us that $\gamma$ extends to an isomorphism 
$\Gamma: \Abar \isom \Mat(n,\Kbar)$. So 
$\sigma(\Gamma) \Gamma^{-1}$ is conjugation by any lift of
$\chi_E(\xi_\sigma)$ to $\GL_n(\Kbar)$.
It follows that $\chi_{E,*}(\xi)$ represents the class of $A$ 
in $H^1(K,\PGL_n)$.
\endproof

The following variant on the above terminology is often helpful.

\begin{Definition}
\label{repdef}
Let $\Theta$ be a theta group for $E[n]$. As a special case
of Definition~\ref{defemb} we define a {\em representation} of $\Theta$
to be an embedding of $\Theta$ in the matrix algebra $\Mat_n(K)$.
In other words, a representation is a morphism of group varieties
$\Theta \to \GL_n$ that preserves scalars.
We recognise diagrams~(\ref{ThetaEdiag}) and~(\ref{Thetadiag}) 
as representations of theta groups.
\end{Definition}

\section{Recovering explicit equations}

Given $\rho \in (R \otimes R)^\times$ representing an 
element $\rho \, \partial R^\times \in H \isom H^1(K,E[n])$ with 
trivial obstruction, we aim to find equations for 
the corresponding Brauer-Severi diagram $[C \to \PP^{n-1}]$.
We present three algorithms for performing this conversion,
assuming in each case the existence of a ``Black Box'' to trivialise 
the obstruction algebra. 
We assume for ease of exposition that $n \ge 3$.

We fix a basis $r_1, \ldots ,r_{n^2}$ for $R$ as a $K$-vector space.
Let $r_1^*, \ldots, r_{n^2}^*$ be the dual basis 
with respect to the trace form
$$ \begin{array}{rcl} 
(r,s) & \mapsto & \tr_{R/K}(rs) = \sum_{T \in E[n]} r(T) s(T).
\end{array} $$
A useful technique, used in several proofs, is reduction to the
``geometric case''. By this we mean taking $K=\Kbar$ and 
$r_i = r_i^* = \delta_{T_i}$ where $E[n] = \{T_1, \ldots ,T_{n^2}\}$.
The following lemma is typical.

\begin{Lemma} 
\label{delta}
Let $\delta = \sum_{i=1}^{n^2} r_i^* \otimes r_i \in R \otimes R$. Then 
$$ \delta(S,T) = \left\{ \begin{array}{ll} 1 & \text{ if } S=T, \\
0 & \text{ otherwise.} \end{array} \right. $$
\end{Lemma} 
\begin{Proof}
We first note that $\delta$ does not depend on the choice of
basis $r_1, \ldots, r_{n^2}$. The lemma follows by reduction to 
the geometric case.
\end{Proof}

\subsection{The Hesse pencil method}
We start with the Hesse pencil method, since it is the simplest
of our three methods both to explain and to implement. 

\begin{Proposition}
\label{hp1}
Let $\Theta$ be the twist of $\Theta_E$ by $\xi \in H^1(K,E[n])$. \\
(i) The theta group $\Theta$ has a representation 
\begin{equation}
\label{thetarep}
 \xymatrix{
0 \ar[r] & \Gm \ar[r] \ar@{=}[d] & \Theta 
\ar[r] \ar[d] & E[n] \ar[r] \ar[d]^{\chi} & 0  \\
0 \ar[r] & \Gm \ar[r]      & \GL_n  \ar[r]    & \PGL_n
\ar[r]           & 0  } 
\end{equation}
if and only if $\Ob(\xi)=0$. \\
(ii) If $\Theta$ has a representation~(\ref{thetarep}) then there is 
a unique genus one normal curve $C \subset \PP^{n-1}$ with Jacobian $E$
for which the action of each $T \in E[n]$ on $C$ is given by $\chi(T)$.
Moreover $[C \to \PP^{n-1}]$ is the Brauer-Severi diagram determined
by $\xi$. 
\end{Proposition}
\begin{Proof}
(i) This is a special case of Proposition~\ref{newob}. \\
(ii) Let $[C \to \PP^{n-1}]$ be the Brauer-Severi diagram determined
by $\xi$. We recall from \S\ref{thetagps} that $\Theta$ has a representation
\begin{equation}
\label{thetarepC}
 \xymatrix{
0 \ar[r] & \Gm \ar[r] \ar@{=}[d] & \Theta 
\ar[r] \ar[d] & E[n] \ar[r] \ar[d]^{\chi_C} & 0  \\
0 \ar[r] & \Gm \ar[r]      & \GL_n  \ar[r]    & \PGL_n
\ar[r]           & 0  } 
\end{equation}
According to Lemma~\ref{natural}, the representations~(\ref{thetarep}) 
and~(\ref{thetarepC}) differ only by an automorphism of $\Mat_n(K)$. 
By the Noether-Skolem theorem this automorphism is 
conjugation by an element of $\GL_n(K)$.
So making a change of co-ordinates on $\PP^{n-1}$ we may 
arrange that $\chi=\chi_C$.

The proof of Lemma~\ref{lem1.8} shows 
that if $C,C' \subset \PP^{n-1}$ are genus one normal curves
with the same $j$-invariant, then there exists $\alpha \in \PGL_n(\Kbar)$ 
with $\alpha(C)=C'$. Moreover, using Lemma~\ref{standardtheta},
one can show that the image of $\chi_C$ is its own centraliser in $\PGL_n$. 
The uniqueness statement follows.
\end{Proof}

We describe the Hesse pencil method in greater detail.

\begin{Proposition}
\label{hp2}
Let $\xi \in H^1(K,E[n])$ and $\rho \in (R \otimes R)^\times$ 
with $w_2(\xi)= \rho \, \partial R^\times$. Let $A_\rho = (R,+,*_{\eps \rho})$
where $\inv_2(\Theta_E)= \eps \, \partial R^\times$. Then \\
(i) $\Ob(\xi) = 0$ if and only if $A_\rho \isom \Mat_n(K)$. \\
(ii) If $\tau : A_\rho \isom \Mat_n(K)$ is an isomorphism of $K$-algebras
and 
$$ \begin{array}{c} M = \sum_{i=1}^{n^2} r_i^* \tau(r_i) \in \GL_n(R) 
= \Map_K(E[n],\GL_n(\Kbar)) \end{array} $$
then there is a unique genus one normal curve $C \subset \PP^{n-1}$ 
with Jacobian $E$ for which the action of each $T \in E[n]$ on $C$ is 
given by $M_T$. Moreover $[C \to \PP^{n-1}]$ is the Brauer-Severi 
diagram determined by $\xi$. 
\end{Proposition}

\begin{Proof}
(i) Let $\Theta$ be the twist of $\Theta_E$ by $\xi$. By Lemma~\ref{thetarho}
we have $\inv_2(\Theta) = \eps \rho \, \partial R^\times$. Then
Lemma~\ref{mrho} identifies $A_\rho$ as the enveloping algebra of
$\Theta$. We are done by Proposition~\ref{newob}. \\
(ii) Since $\inv_2(\Theta) = \eps \rho \, \partial R^\times$ there exists
a Galois equivariant section $\phi : E[n] \to \Theta$ with
$$\phi(S) \phi(T) = \eps(S,T) \rho(S,T) \phi(S+T) $$
for all $S,T \in E[n]$. We claim that $\phi(T) \mapsto M_T$ extends
to a representation of $\Theta$. 
It suffices to check this in the geometric case, 
whereupon $M_T = \tau(\delta_T)$. 
The proof of Lemma~\ref{mrho} shows that $\phi(T) \mapsto \delta_T$ 
extends to an embedding of $\Theta$ in $A_\rho$. Since $\tau$
is an isomorphism of $K$-algebras, it follows that $\Theta$ 
embeds in $\Mat_n(K)$ as claimed.

Finally we apply Proposition~\ref{hp1}(ii) with $\chi(T)= [M_T]$. 
\end{Proof}

\begin{Remark} We may take any convenient choice for $\eps$. For
example we saw in Lemma~\ref{thetaodd} that if $n$ is odd then 
a convenient choice is the square root of the Weil pairing.
\end{Remark}

It remains to recover equations for $C \subset \PP^{n-1}$ 
from $M \in \GL_n(R)$. 

\begin{Proposition}
\label{Yn}
Assume $n \ge 3$ and let $M \in \GL_n(R)$ such that \\
(i) $E[n] \to \PGL_n \, ; \, T \mapsto [M_T]$ 
is a group homomorphism, and \\
(ii) $M_S M_T M_S^{-1} M_T^{-1} = e_n(S,T) I_n$ for all $S,T \in E[n]$. \\
Then the genus one normal curves $C \subset \PP^{n-1}$ for which each
matrix $M_T$ acts as translation by some $n$-torsion point of $\Jac(C)$,
are parametrised by (a twist of) the modular curve $Y(n)$. Moreover
the number of curves in this family that are defined over $K$ 
and have Jacobian $E$ is 
$$ \nu_{E,n} = [\Aut_K(E[n]): \Aut_K(E)] $$
where $\Aut_K(E[n])$ is the group of $K$-rational
automorphisms of $E[n]$ that respect the Weil pairing. 
\end{Proposition}
\begin{Proof}
The first statement is a geometric one. 
For its proof we may fix a basis $S,T$ for $E[n]$ and assume 
that $M_S$ and $M_T$ are the standard matrices $M_1$ and $M_2$ 
specified in Lemma~\ref{standardtheta}.
We recall that $Y(n)$ parametrises the triples $(E',S',T')$ where
$E'$ is an elliptic curve and $S',T'$ are a basis for $E'[n]$ with
$e_n(S',T') = \zeta_n$. Lemma~\ref{standardtheta} furnishes 
us with a bijection between the triples $(E',S',T')$ and the genus 
one normal curves $C \subset \PP^{n-1}$ considered here. 
So the latter are also parametrised by $Y(n)$. 

Proposition~\ref{hp1}(ii) establishes the existence of a genus one normal
curve $C \subset \PP^{n-1}$ with Jacobian $E$ for which
the action of each $T \in E[n]$ is given by $M_T$. Let $C'$ be another
curve in the $Y(n)$ family, defined over $K$ and with 
$\Jac(C') \isom E$. Then each $T \in E[n]$ acts on $C'$ via 
$M_{\alpha(T)}$ for some $\alpha \in \Aut_K(E[n])$. Changing 
our choice of isomorphism $\Jac(C') \isom E$ changes $\alpha$ by
an element of $\Aut_K(E)$. Here we use our assumption $n \ge 3$ to
identify $\Aut_K(E)$ as a subgroup of $\Aut_K(E[n])$.
It follows by Proposition~\ref{hp1}(ii) that the curves $C'$
considered here are in bijection with the quotient 
group $ \Aut_K(E[n])/\Aut_K(E)$.
\end{Proof}

In general we have $\nu_{E,n} \le \# \PSL_2(\Z/n\Z) $.
We now specialise to the case $n=3$.

\begin{Lemma}
Let $M  \in \GL_3(R)$ as in Proposition~\ref{Yn}.  
Let $(x:y:z)$ be co-ordinates on $\PP^2$ and put
$$ \begin{pmatrix} x' \\ y' \\ z' \end{pmatrix} = M 
 \begin{pmatrix} x \\ y \\ z \end{pmatrix}, \qquad 
\begin{pmatrix} x'' \\ y'' \\ z'' \end{pmatrix} = M^2 
 \begin{pmatrix} x \\ y \\ z \end{pmatrix}. $$
Then the family of curves parametrised by $Y(3)$ is (an open subset of) 
the pencil of plane cubics spanned by the $F_i(x,y,z) \in K[x,y,z]$ where
$$ \left| \begin{array}{ccc} x & x' & x'' \\ y & y' & y'' \\ z & z' & z''
\end{array} \right|  = \sum_{i=1}^9 F_i(x,y,z) r_i. $$
\end{Lemma}
\begin{Proof} 
Let $S,T$ be a basis for $E[n]$. For the proof we may
assume that $M_S$ and $M_T$ are the standard matrices $M_1$ and $M_2$ 
specified in Lemma~\ref{standardtheta}. In this case our construction 
does indeed give the pencil of plane cubics 
spanned by $x^3+y^3+z^3$ and $xyz$. 
\end{Proof}

We use the classical invariants of a ternary cubic
(cf. \cite{AKM3P}) to pick out those members of the pencil that
are defined over $K$ and have Jacobian $E$. 
In practice this means finding the $K$-rational roots of a polynomial
of degree 12. According to 
Proposition~\ref{Yn} we are left with a list of $\nu_{E,3}$ 
candidates for $C \subset \PP^2$. In favourable circumstances 
(including the case $K=\Q$) we can show $\nu_{E,3} = 1$.
\begin{Lemma} Assume that either  
(i) $\rho_{E,3} : G_K \to \GL_2(\Z/3\Z)$ is surjective, or 
(ii) $K$ is a number field with a real place. 
Then $\nu_{E,3} = 1$.
\end{Lemma}
\begin{Proof}
Let $G \subset \GL_2(\Z/3\Z)$ be the image of $\rho_{E,3}$. Then
\begin{equation}
\label{cent}
 \Aut_K(E[3]) = \{ x \in \SL_2(\Z/3\Z) \mid xy =yx \text{ for all }
y \in G \}.
\end{equation}
(i) We find  $\Aut_K(E[3]) = \{\pm 1\}$ and so $\nu_{E,3}= 1$. \\
(ii) Since $\zeta_3 \not\in K$ it is clear that $G$ contains an 
element of determinant $-1$. But there are only $3$ such conjugacy classes 
in $\GL_2(\Z/3\Z)$. Hence we may assume that
$G$ contains either $\pm a$ or $b$ where 
$$ a = \begin{pmatrix} 1 & 1 \\ -1 & 1 \end{pmatrix}
\quad \text{ and } \quad
b = \begin{pmatrix} 1 & 0 \\ 0 & -1 \end{pmatrix}. $$
It follows by~(\ref{cent}) that 
$\Aut_K(E[3]) = \{\pm 1\}$ or $\{\pm 1,\pm a^2\}$. If the latter,
then a further application of~(\ref{cent}) shows that $G$ is cyclic of
order $8$. By considering the subfields of $K(E[3])$ we are led 
to the contradiction $K(\sqrt{-3}) \subset K(E[3]) \cap \R$. 
Hence $\Aut_K(E[3]) = \{\pm 1\}$ and $\nu_{E,3} = 1$ as claimed.
\end{Proof}

In \cite{testeq} a formula based on Corollary~\ref{detMcor} is used to
recover $\alpha \in R^\times/(R^\times)^3$ from a ternary cubic.
The Hesse pencil method is completed in the case $\nu_{E,3} > 1$
by applying this formula to each of the (at most 12) candidate 
ternary cubics, and seeing which gives rise to the correct element 
$\alpha \in R^\times/(R^\times)^3$.

\subsection{The flex algebra method} 
This method has the advantage over the Hesse pencil method that it
works for all $n \ge 2$. For ease of exposition we continue to
assume that $n \ge 3$.

We embed $E \to \PP^{n-1}$ via the complete linear system $|n.0|$
and compute $M \in \GL_n(R) = \Map_K(E[n],\GL_n(\Kbar))$ 
describing the action of $E[n]$ on $E$ by translation.
We also compute $\eps \in (R \otimes R)^\times$ with $M_S M_T = \eps(S,T)
M_{S+T}$ for all $S,T \in E[n]$.

\begin{Proposition}
\label{fa1}
Let $\xi \in H^1(K,E[n])$ and $\rho \in (R \otimes R)^\times$ 
with $w_2(\xi)= \rho \, \partial R^\times$. 
Let $A_1 = (R,+,*_{\eps})$ and $A_\rho = (R,+,*_{\eps \rho})$. \\
(i) There exists $\gamma \in \Rbar^\times$ with $\partial \gamma=\rho$,
and so an isomorphism of $\Kbar$-algebras $.\gamma : \Abar_\rho \to
\Abar_1$. \\
(ii) The map $\tau_1 : A_1 \to \Mat_n(K)$ given by
$\tau_1(x)_{ij} = \tr_{R/K} (x M_{ij})$ is an isomorphism of $K$-algebras. \\
(iii) If $\tau_\rho : A_\rho \isom \Mat_n(K)$ is an isomorphism of $K$-algebras
then there is a commutative diagram
\begin{equation}
\label{fa1diag}
 \xymatrix{ \Abar_\rho \ar[r]^-{\tau_\rho} \ar[d]_{.\gamma} & \Mat_n(\Kbar) 
 \ar[d]^{\beta} \\  \Abar_1  \ar[r]^-{\tau_1} & \Mat_n(\Kbar) } 
\end{equation}
where $\beta$ is conjugation by some matrix $B \in \GL_n(\Kbar)$.
Moreover $B$ represents a change of co-ordinates on $\PP^{n-1}$
taking the Brauer-Severi diagram $[E \to \PP^{n-1}]$ to its twist 
$[C \to \PP^{n-1}]$ by $\xi$.
\end{Proposition}

\begin{Proof}
(i) The element $\gamma \in \Rbar^\times$ exists by Corollary~\ref{kbarcor}.
The isomorphism $.\gamma$ is that specified in Lemma~\ref{concrete}. \\
(ii) We must show that $\tau_1$ is a ring homomorphism. 
Reducing to the geometric case we have $\tau_1(\delta_T)=M_T$. 
Since $\delta_S *_\eps \delta_T = \eps(S,T) \delta_{S+T}$ 
and $M_S M_T = \eps(S,T) M_{S+T}$ the result is clear. \\
(iii) Let $\beta$ be the isomorphism of $\Kbar$-algebras 
making~(\ref{fa1diag}) commute. By the Noether-Skolem theorem it is
conjugation by some matrix $B \in \GL_n(\Kbar)$.

Let $\Theta_E$ and $\Theta$ be the theta groups for $E[n]$ with 
second invariants $\inv_2(\Theta_E) = \eps \partial R^\times$
and $\inv_2(\Theta) = \eps \rho \partial R^\times$. 
By Lemma~\ref{mrho} the enveloping algebras are $A_1$ and $A_\rho$.
We interpret the isomorphisms $\tau_1: A_1 \isom \Mat_n(K)$ 
and $\tau_{\rho} : A_\rho \isom \Mat_n(K)$ as representations of 
$\Theta_E$ and $\Theta$. 
So $\Theta_E$ and $\Theta$ are now subgroups
of $\GL_n$ generated up to scalars by the $\tau_1(\delta_T)$,
respectively $\tau_\rho(\delta_T)$, for $T \in E[n]$. 

Since $\tau_1(\delta_T)= M_T$ the theta group $\Theta_E \subset \GL_n$
is that determined by $[E \to \PP^{n-1}]$.
On the other hand Proposition~\ref{hp1} tells us that $\Theta \subset \GL_n$
is the theta group for some $[C \to \PP^{n-1}]$. 
Moreover, since $\Theta$ is the twist of $\Theta_E$ by $\xi$, the 
proposition also tells us that $[C \to \PP^{n-1}]$ is the twist of 
$[E \to \PP^{n-1}]$ by $\xi$. The commutativity of~(\ref{fa1diag}) 
shows that $\Theta_E \subset \GL_n$ and $\Theta \subset \GL_n$ 
are related by conjugation by $B$. It follows by the uniqueness 
statement of Proposition~\ref{hp1} that $B$ represents a change 
of co-ordinates on $\PP^{n-1}$ taking $[C \to \PP^{n-1}]$ 
to $[E \to \PP^{n-1}]$.
\end{Proof}

As it stands the method is unsatisfactory, since we have to solve
for $\gamma \in \Rbar^\times = (R \otimes \Kbar)^\times$ 
with $\partial \gamma = \rho$.
By~(\ref{exseq1}) the $n^2$ choices for $\gamma$ form a coset of $E[n]$
inside $\Rbar^\times$. This coset is an $E[n]$-torsor,
which turns out to be the twist of $E[n]$ by $\xi$. 
A glance at Remark~\ref{remflex} now suggests we should solve
for $\gamma \in (R \otimes F)^\times$ where $F$ is the field
of definition of a flex point on $C$. This is already clear
from the conclusions of Proposition~\ref{fa1}, since the $K$-rational
point $0 \in E$ gets mapped to a flex point on $C$.
These observations motivate the following 
refinement of Proposition~\ref{fa1}.

\begin{Proposition}
\label{fa2}
Let $\xi \in H^1(K,E[n])$ and $\rho \in (R \otimes R)^\times$ 
with $w_2(\xi)= \rho \, \partial R^\times$. 
Let $A_1 = (R,+,*_{\eps})$, $A_\rho = (R,+,*_{\eps \rho})$
and $F = (R,+,*_{\rho})$. \\
(i) Let $\Phi$ be the $E[n]$-torsor determined by $\xi$.
Then $F$ is the \'etale algebra of $\Phi$.
In particular $F$ is a product of field extensions of $K$. \\
(ii) There is an isomorphism of $F$-algebras
$$ \begin{array}{lrcl}
\alpha : & A_\rho \otimes_K F & \to & A_1 \otimes_K F \\
& x \otimes 1 & \mapsto & \sum_{i=1}^{n^2} 
 r^*_i x  \otimes r_i. 
\end{array} $$
(iii) Let $\tau_1 : A_1 \isom \Mat_n(K)$ and $\tau_\rho: A_\rho
\isom \Mat_n(K)$ be the isomorphisms of Proposition~\ref{fa1}. Then
there is a commutative diagram
\begin{equation}
\label{alphabeta}
 \xymatrix{ A_\rho \otimes_K F \ar[d]_\alpha 
\ar[r]^{\tau_\rho}  & \Mat_n(F)  \ar[d]^\beta \\ 
A_1 \otimes_K F \ar[r]^{\tau_1} & \Mat_n(F) }
\end{equation} 
where $\beta$ is conjugation by some matrix $B \in \GL_n(F) 
=\Map_K(\Phi,\GL_n(\Kbar))$. 
Moreover for each $P \in \Phi$ the matrix $B_P \in \GL_n(\Kbar)$ 
represents a change of co-ordinates on $\PP^{n-1}$
taking the Brauer-Severi diagram $[E \to \PP^{n-1}]$ 
to its twist $[C \to \PP^{n-1}]$ by $\xi$.
\end{Proposition}
\begin{Proof} 
(i) This is proved in Proposition~\ref{flexenv}. \\
(ii) We first show that $\alpha$ is a ring homomorphism, i.e.
\begin{equation}
\label{homcheck}
\begin{array}{rcl}
\sum_{i=1}^{n^2} r^*_i(x *_{\eps \rho} y)
\otimes r_i 
&  = &  \sum_{i,j=1}^{n^2} 
( r^*_i x *_\eps r^*_j y) \otimes ( r_i *_\rho r_j) 
\end{array} 
\end{equation}
for all $x,y \in R$. Since $\alpha$ does not depend on the choice
of basis $r_1, \ldots, r_{n^2}$ we may reduce to the geometric case.
Putting $x = \delta_S$ and $y = \delta_T$ in~(\ref{homcheck}) it becomes 
$$ \begin{array}{rcl}
  \eps(S,T) \rho(S,T) \delta_{S+T}  \otimes 
 \delta_{S+T}   &  = &  \eps(S,T) 
\delta_{S+T}  \otimes \rho(S,T) \delta_{S+T} 
\end{array} $$
which is a tautology. 
Since $A_\rho$ and $A_1$ are central simple algebras (of the same dimension) 
it follows that $\alpha$ is an isomorphism. \\
(iii) Let $\beta$ be the isomorphism of $F$-algebras making the
diagram commute. By the Noether-Skolem theorem (applied to each
constituent field of $F$) it is conjugation by some matrix $B \in \GL_n(F)$.

Let $\gamma = ( 1 \otimes \iota_F) (\delta) \in R \otimes F$ where
$\delta \in R \otimes R$ was defined in Lemma~\ref{delta},
and $\iota_F : R \isom F$ is the isomorphism of underlying $K$-vector spaces. 
Then the isomorphism $\alpha : A_\rho \otimes F \isom A_1 \otimes F$
is multiplication by $\gamma$ in $R \otimes F$. So for each $P \in \Phi$
we obtain a diagram
\begin{equation*}
 \xymatrix{ \Abar_\rho \ar[r]^-{\tau_\rho} \ar[d]_{.\gamma_P} & \Mat_n(\Kbar) 
 \ar[d]^{\beta_P} \\  \Abar_1  \ar[r]^-{\tau_1} & \Mat_n(\Kbar) } 
\end{equation*}
where $.\gamma_P$ is multiplication by $\gamma_P$ in $\Rbar^\times$
and $\beta_P$ is conjugation by $B_P \in \GL_n(\Kbar)$.
We are done by Proposition~\ref{fa1}.
\end{Proof}

We summarise the flex algebra method in the following 7 steps.

\medskip

\paragraph{Step 1}
We embed $E \to \PP^{n-1}$ via the complete linear system $|n.0|$. 
The curve $E$ is now defined by homogeneous polynomials $f_1, \ldots, f_N$
in $K[x_1, \ldots,x_n]$.

\medskip

\paragraph{Step 2}
We compute $M \in \GL_n(R)$ and $\eps \in (R \otimes R)^\times$ 
as described at the start of this subsection.
If $n$ is odd then, in the notation of Lemma~\ref{thetaodd}, 
we may choose $M$ with $M^n=I_n$ and 
$\iota M \iota^{-1} = M^{-1}$. This enables us to take 
$\eps = e^{1/2}$.

\medskip

\paragraph{Step 3}
Let $A_1= (R,+,*_\eps)$. We compute the isomorphism of $K$-algebras 
$\tau_1 : A_1 \isom \Mat_n(K)$ specified in Proposition~\ref{fa1}.

\medskip

\paragraph{Step 4}
Let $A_\rho = (R,+,*_{\eps \rho})$. 
We use the
Black Box to find an isomorphism of $K$-algebras $\tau_{\rho} : A_\rho
\isom \Mat_n(K)$. 

\medskip

\paragraph{Step 5}
Let $F = (R,+,*_\rho)$.
We compute the composite
$$ \tau'_\rho : A_\rho \stackrel{\alpha}{\ra} 
A_1 \otimes F \stackrel{\tau_1}{\ra} \Mat_n(K) \otimes F = \Mat_n(F). $$
It is given by 
$\tau_{\rho}'(x) = \sum_{i=1}^{n^2} \tau_1 (r^*_i x) \otimes r_i.$

\medskip

\paragraph{Step 6}
We use linear algebra to solve for $B \in \GL_n(F)$ with 
$\tau'_\rho(x) = B \, \tau_\rho(x) \, B^{-1}$ for all $x \in R$. 

\medskip

\paragraph{Step 7}
Let $w_1, \ldots ,w_{n^2}$ be a basis for $F$ as a $K$-vector
space. Then $C$ is defined by the homogeneous polynomials
$g_{ij} \in K[x_1, \ldots ,x_n]$ with
$$ \begin{array}{rcl}
f_i(\sum_{j=1}^n B_{1j} x_j, \ldots , \sum_{j=1}^n B_{nj} x_j )
 & = & \sum_{j=1}^{n^2} w_j g_{ij} (x_1, \ldots,x_n). 
\end{array} $$

\begin{Remark} In Steps 6 and 7 it suffices to work with
any constituent field of $F$. But in the generic case 
Galois acts transitively on the flex points of $C$. So $F$ is already
a field and there is no saving to be made.
\end{Remark}

\subsection{The Segre embedding method}
The third of our algorithms leads more directly to equations for
$C$ (and avoids the need to compute the flex algebra). 
Here we confine ourselves to a brief description.  
Further details, including a proof that the method works, 
will be given in the second paper of this series \cite{paperII}.

Let $\rho \, \partial R^\times \in H$ correspond to a torsor divisor
class pair $(C,[D])$. Even before we use the Black Box, we can 
write down equations for $C$ as a genus one normal curve in 
$\PP^{n^2-1}$ with hyperplane section $nD$. 
To do this we fix a Weierstrass equation for $E$ and write  $(x(P),y(P))$ 
for the co-ordinates of $P \in E \setminus \{0\}$. We also write 
$\lambda(P_1,P_2)$ for the slope of the chord through $P_1, P_2  
\in E \setminus \{0\}$ with $P_1+P_2 \not= 0$,
respectively of the tangent line if $P_1 = P_2$. 
We put $z = \sum r_i z_i$ where  $z_1, \ldots , z_{n^2}$ are indeterminates.
Since $R = \Map_K(E[n],\Kbar)$ we have 
$$ \begin{array}{c}  z(T)  =  
\sum_{i=1}^{n^2} r_i(T) z_i \in \Kbar[z_1, \ldots , z_{n^2}]. 
\end{array} $$

For $T \in E[n] \setminus \{0\}$ we consider the polynomial 
$$ \big(x-x(T)\big)z(0)^2 - \rho(T,-T) z(T) z(-T) $$ 
in $\Kbar[x,z_1, \ldots,z_{n^2}]$. 
We define a quadric of type 1 to be the difference of any two 
such polynomials. These quadrics span a $\Kbar$-vector subspace 
of $\Kbar[z_1, \ldots, z_{n^2}]$ of dimension $d_1$ where 
$$ d_1 = \left\{ \begin{array}{ll} (n^2-3)/2 & \text{ if $n$ odd,} \\
n^2/2 & \text { if $n$ even.} \end{array} \right. $$

For $T,T_1,T_2 \in E[n] \setminus \{0\}$ with 
$T_1+T_2=T$ we consider the polynomial
$$ \big(\lambda_T-\lambda(T_1,T_2)\big)z(0)z(T) 
- \rho(T_1,T_2) z(T_1) z(T_2) $$
in $\Kbar[\lambda_T,z_1, \ldots,z_{n^2}]$. 
We define a quadric of type 2 to be the difference of any two 
such polynomials that share the same choice of $T$. These quadrics 
span a $\Kbar$-vector subspace 
of $\Kbar[z_1, \ldots, z_{n^2}]$ of dimension $d_2$ where 
$$ d_2 = \left\{ \begin{array}{ll} (n^2-1)(n^2-3)/2 & \text{ if $n$ odd,} \\
n^2(n^2-4)/2 & \text { if $n$ even.} \end{array} \right. $$

It is clear that the spaces of quadrics of types 1 and 2 are each 
Galois invariant. We thus obtain a $K$-vector space of quadrics in 
$K[z_1, \ldots ,z_{n^2}]$ of dimension $d_1+d_2 = n^2(n^2-3)/2$. 
It is shown in the second paper of this series \cite{paperII} 
that these quadrics generate the homogeneous ideal 
of $C \subset \PP^{n^2-1}$ embedded by $|nD|$. Notice that 
our formulae for the dimension of each vector space 
of quadrics are easily checked by reduction to the geometric case.

\begin{Remarks}
(i) In practice when computing the quadrics defining $C$ we work 
over the constituent fields of $R \otimes R$, rather than over $\Kbar$. \\ 
(ii) It is arguably more natural first to give equations for 
$C \subset \PP(R)$ and only then to identify $\PP(R) = \PP^{n^2-1}$
by means of our choice of basis $r_1, \ldots , r_{n^2}$ for $R$.  \\
(iii) We recall that $\rho \, \partial R^\times$ is an element of $H \subset
(R \otimes R)^\times / \partial R^\times$. If we multiply $\rho$
by $\partial \gamma$ for some $\gamma \in R^\times$ then the
effect on $C$ is that of a change of co-ordinates on $\PP^{n^2-1}$. \\
(iv) In the case $n=2$ we find that $C \subset \PP^3$ is the
complete intersection of two quadrics of type 1. Our method 
reduces to the classical number field method for 2-descent.
\end{Remarks}

We now have equations for $C \subset \PP^{n^2-1}$ embedded by $|nD|$.
It remains to compute equations for $C \subset \PP^{n-1}$ 
embedded by $|D|$. This is achieved in the following 5 steps.

\medskip

\paragraph{Step 1}
We compute $F \in R(E)^\times = \Map_K(E[n],\Kbar(E)^\times)$ 
with $\divv(F_T) = n.T-n.0$. We fix a local parameter $z$ in the
local ring of $E$ at $0$ and scale each of the rational functions $F_T$ 
to have leading coefficient 1 when expanded as a Laurent power series in $z$. 

\medskip

\paragraph{Step 2}
We compute $\eps \in (R \otimes R)^\times$ with 
$$ \eps(T_1,T_2) = \frac{ F_{T_1+T_2} (P)}{F_{T_1}(P)F_{T_2}(P-T_1)}$$
for all $T_1,T_2 \in E[n]$ and $P \in E \setminus \{0,T_1,T_1+T_2\}$. 
(This formula for $\eps$ was derived near the end of \S\ref{etal}.)

\medskip

\paragraph{Step 3}
Let $A_\rho = (R,+,*_{\eps \rho})$.
We use the Black Box to find an isomorphism of $K$-algebras $\tau : A_\rho
\isom \Mat_n(K)$.

\medskip

\paragraph{Step 4}
For each quadric $f(z_1, \ldots, z_{n^2})$ computed above
we make a change of co-ordinates to obtain a quadric 
$g(x_{11},x_{12}, \ldots,x_{nn})$ with
$$ \begin{array}{rcl} 
g ( \sum_{i=1}^{n^2} \tau(r_i)_{11} z_i , \ldots, 
\sum_{i=1}^{n^2} \tau(r_i)_{nn} z_i) & = & f(z_1, \ldots,z_{n^2}). 
\end{array} $$
The new quadrics define $C$ as a 
curve 
in $\PP(\Mat_n)$.

\medskip

\paragraph{Step 5}
There is a direct sum decomposition $\Mat_n = \langle I_n \rangle 
\oplus \{ \Tr = 0 \}$ where $\{ \Tr=0 \}$ is the subspace of matrices
of trace zero. We project $C$ to the trace zero subspace.
Then $C$ is contained in the rank 1 locus. In other words $C$ lies
in the image of the Segre embedding
$$ \PP^{n-1} \times (\PP^{n-1})^\vee \to \PP(\Mat_n). $$
We pull back to a curve in $\PP^{n-1} \times (\PP^{n-1})^\vee$
and finally project onto the first factor to obtain 
$C \to \PP^{n-1}$. (In fact, projecting onto the second 
factor gives the dual curve.)

\begin{Remark} In both the flex algebra method and the Segre
embedding method we have specified particular choices
for $\eps \in (R \otimes R)^\times$. But for the purposes
of computing the obstruction algebra, we can use any
$\eps \in (R \otimes R)^\times$ with 
$\inv_2(\Theta_E) = \eps \, \partial R^\times$.
Then after trivialising the obstruction algebra
we can use the isomorphism of Lemma~\ref{concrete}
to compensate for having made a different choice of $\eps$.
So the method we use (Hesse pencil, flex algebra
or Segre embedding) has no impact on the implementation
of the Black Box.
\end{Remark}

\end{document}